\documentclass[final,leqno,onefignum,onetabnum]{siamltex1213RBP}

\usepackage{benstyle_SIAM}
\newcommand*{\asinh}{\mathrm{asinh}}

\numberwithin{equation}{section}
\numberwithin{theorem}{section}

\usepackage{overpic}
\def\complex{{\mathbb C}}

\def\psiDE{\psi_{DE}}		
\def\psiDEi{\psi_{DE}^{-1}}		
		
\def\psiSE{\psi_{SE}^{}}		
\def\psiSEi{\psi_{SE}^{-1}}		
\def\psiSDE{\psi_{SDE}^{}}		
\def\psiSDEi{\psi_{SDE}^{-1}}

\def\s{\hspace{-1pt}}		
		
\def\ee{\hspace{-2pt}=\hspace{-2pt}}

\title{Resolution-optimal exponential and double-exponential transform methods for functions with endpoint singularities}

\author{Ben Adcock%
        \thanks{Department of Mathematics,	
		Simon Fraser University,
		8888 University Drive,
		Burnaby,
		BC V5A 1S6,
		Canada
                  (\texttt{ben\_adcock@sfu.ca},
                   \texttt{http://www.sfu.ca/\textasciitilde benadcock/})
                   }
                   \and
              Jes\'us Mart\'in--Vaquero%
                   \thanks{
                    Department of Applied Mathematics, E.T.S.I.I. B\'ejar,
 University of Salamanca, E37700 Spain
   (\texttt{jesmarva@usal.es},
                   \texttt{http://diarium.usal.es/jesmarva/})
                   }
                   \and
             Mark Richardson%
                   \thanks{
                   21 Twyford Court,
		  Fortis Green,
                   London N10 3ES,
                   United Kingdom
                   (\texttt{mrichardson82@gmail.com})
                   }
                   }

\begin{document}

\maketitle
\slugger{sisc}{xxxx}{xx}{x}{x--x}

\begin{abstract}
We introduce a numerical method for the approximation of functions which are analytic on compact intervals, except at the endpoints.  This method is based on variable transforms using particular parametrized exponential and double-exponential mappings, in combination with Fourier-like approximation in a truncated domain.  We show theoretically that this method is superior to variable transform techniques based on the standard exponential and double-exponential mappings.  In particular, it can resolve oscillatory behaviour using near-optimal degrees of freedom, whereas the standard mappings require degrees of freedom that grow superlinearly with the frequency of oscillation.  We highlight these results with several numerical experiments.  Therein it is observed that near-machine epsilon accuracy is achieved using a number of degrees of freedom that is between four and ten times smaller than those of existing techniques.
\end{abstract}

\begin{keywords}
endpoint singularities, conformal maps, Fourier series, resolution analysis
\end{keywords}

\begin{AMS}
41A05, 41A25, 65D05
\end{AMS}

\pagestyle{myheadings}
\thispagestyle{plain}
\markboth{B.\ ADCOCK, J.\ MART\'IN--VAQUERO AND M.\ RICHARDSON}{RESOLUTION-OPTIMAL TRANSFORM METHODS}

\section{Introduction}
Analytic functions on compact intervals can be accurately approximated using either Fourier (in the case of periodic functions) or Chebyshev expansions.  Both approximations can be computed efficiently via the Fast Fourier Transform (FFT) and are known to converge geometrically fast in the number of degrees of freedom \cite{boyd,TrefethenATAP}.  For periodic functions, Fourier expansions are typically preferable over Chebyshev expansions, since they are more efficient by a factor of $\pi/2$ at resolving oscillatiory behaviour.

In this paper, we consider the fast and accurate approximation of functions which are analytic on compact intervals, except possibly at the endpoints.  Such functions arise in a number of applications in scientific computing, and their accurate approximation via so-called \textit{variable transform methods} has been the subject of a long line of research (see \cite{LundBowers,StengerSiam81,StengerBook93,StengerBook10} and references therein).  These methods are typically based on the following approach.  First, a function $f(x)$ on a compact interval $[0,1]$ (without loss of generality) is transformed via an invertible mapping $\psi : (0,1) \rightarrow (-\infty,\infty)$ to a function $F(s) = f(\psi^{-1}(s))$ defined over the real line.  Subject to mild smoothness assumptions on $f$, standard choices of $\psi$, based on \textit{exponential} or \textit{double-exponential} transforms, result in functions $F(s)$ which decay exponentially or double-exponentially fast to their limiting values as $|s| \rightarrow \infty$.  Hence, $F$ can be approximated by \textit{domain truncation}: a parameter $L > 0$ is fixed, and then $F$ is approximated on the interval $[-L,L]$ by a standard technique, e.g.\ Chebyshev expansions.  Provided $L$ is chosen sufficiently large, one can expect a good approximation to $F$, and therefore $f$.  Note that $\mathrm{sinc}$ interpolation is also commonly used in place of domain truncation (see Remark \ref{r:whynotsinc} for a discussion).

In a recent paper by two of the authors, it was shown that the standard exponential mapping $\psi_E$ used in practice has poor resolution properties for oscillatory functions \cite{AdcockRichardsonMappings}.  This analysis was based on domain truncation with Chebyshev interpolation in the truncated interval $[-L,L]$.  It was proved that the number of degrees of freedom required to resolve oscillations of frequency $\omega$ scales not linearly, but quadratically in $\omega$.  This scaling is substantially worse than the case of either Chebyshev or Fourier approximations of analytic functions, both of which are linear with small constants.  Aiming to improve this behaviour a new, parametrized mapping $\psi_{SE}$ was introduced in \cite{AdcockRichardsonMappings}.  When combined with Chebyshev approximation in the truncated interval, it was proved that the resulting approximation was able to achieve not only a linear scaling with $\omega$, but also a resolution constant that was close to that of standard Chebyshev expansions for analytic functions.

The purpose of this paper is twofold.  First, as noted by Boyd \cite{BoydChebyshevInferior}, Chebyshev domain truncation is typically inferior to Fourier domain truncation for approximating analytic functions on infinite intervals.  Hence, seeking to improve the technique of \cite{AdcockRichardsonMappings} even further, we introduce a new approximation strategy for the truncated domain which, similar to the improvement of Fourier over Chebyshev expansions for analytic functions, enhances the resolution power by a factor of $\pi/2$ over the Chebyshev-based approximations considered therein.  This strategy is based on a Fourier-type expansion, and can be implemented efficiently using FFTs.  We show that the resulting numerical methods converge root-exponentially fast in the number of degrees of freedom for both $\psi_{E}$ and $\psi_{SE}$.  However, unlike $\psi_{E}$, the parametrized exponential mapping $\psi_{SE}$ possesses near-optimal resolution power.  As we establish, careful selection of the various parameters allows one to make the resolution constant arbitrarily close to that of Fourier expansions for analytic and periodic functions.

Second, aiming to enhance convergence of these method, we introduce a new parametrized double exponential mapping $\psi_{SDE}$.  Using the same approximation strategy in the truncated domain, we show that this new mapping achieves near-optimal resolution power, much like that of $\psi_{SE}$, but with a convergence rate that is nearly exponential.  This order of convergence is similar to that obtained from the standard double-exponential mapping $\psi_{DE}$, but the resolution power is substantially enhanced.  In particular, it is linear in the frequency $\omega$ with a constant that can be made arbitrarily close to optimal, in comparison to $\ord{\omega \log \omega}$ which, as we show, is the corresponding resolution power for $\psi_{DE}$.

A summary of the convergence and resolution power of the numerical methods introduced in this paper is given in Table \ref{table:resultsSummary}.  We note also that all the methods are fast, and can be implemented in $\ord{n \log n}$ time using FFTs, where $n$ is the number of degrees of freedom.  Moreover, the improved resolution properties of the new mappings $\psi_{SE}$ and $\psi_{SDE}$ lead to significant gains in efficiency.  In our numerical experiments, we present examples where near-machine epsilon accuracy is achieved with these new mappings using a factor of $4-10$ fewer degrees of freedom than when using the standard mappings $\psi_{E}$ and $\psi_{DE}$.

\begin{table}
\centering
\renewcommand{\arraystretch}{1.3} 
\begin{tabular}{c | c c c c}
& $\psi_E$ & $\psi_{SE}$  & $\psi_{DE}$ & $\psi_{SDE}$
\\ \hline
$L$  & $c \sqrt{n}$ & $L_0 + 1/2$ & $1+W(cn)$ & $L_0 + 1/2$ \\
$\alpha$ & --- & $\alpha_0 / \sqrt{n} $& --- & $L_0 \pi / (\pi/2+W(cn))$  \\
error & $\cO_a\left(\rho^{-\sqrt{n}}\right )$ & $\cO_a\left(\rho^{-\sqrt{n}}\right )$ & $\cO_a \left ( \rho^{-n / \log(cn)} \right )$ & $\cO_a \left ( \rho^{-n / \log(cn)} \right )$   \\
d.o.f.& $\ord{\omega^2}$ & $\ord{\omega}$ & $\ord{\omega \log(c\omega)}$ & $\ord{\omega}$ \\
$r$ & $\infty$ & $4 L_0 + 2$ & $\infty$ & $4 L_0 + 2$ \\
\end{tabular}
\caption{\label{table:resultsSummary} \small
A summary of the main results of the paper.  The rows indicate: the parameters $L$ and $\alpha$ (where applicable), the asymptotic approximation error, the number of degrees of freedom required to resolve an oscillation of frequency $\omega$, the asymptotic resolution constant (points-per-wavelength).  Note that $c,\alpha_0,L_0 > 0$  are constants which much be selected by the user.  The function $W$ is the Lambert-W function (see \S \ref{s:DE}).  Our main result is that the new maps $\psi_{SE}$ and $\psi_{SDE}$  achieve the same asymptotic orders of convergence as the standard maps, yet their resolution power is substantially improved.  In particular, the resolution constant $r$ can be made arbitrarily close to $2$, which is the optimal value for resolving oscillatory functions.
}
\end{table}

\subsection*{Outline}  The outline of the remainder of this paper is as follows.  In \S \ref{s:prelims} we introduce variable transform techniques based on domain truncation, the various mappings we consider in this paper and the notion of resolution power.  We introduce the new approximation technique in \S \ref{s:FourierLike} and present a general error analysis.  The next four sections, \S \ref{s:psiE}--\ref{s:psiSDE}, are devoted to proving the main results summarized in Table \ref{table:resultsSummary}.  Finally, in \S \ref{s:numerical} we present numerical experiments.

Throughout this paper, we use the following notation.  We write $B(n) = \cO_a(A(n))$ as $n \rightarrow \infty$ if there exists a $p > 0$ such that $B(n) = \ord{n^p A(n)}$ as $n \rightarrow \infty$.

\section{Preliminaries}\label{s:prelims}
As in \cite {AdcockRichardsonMappings}, let $f(x)$ be a function that is analytic on $(0,1)$ and continuous on $[0,1]$ and suppose that $\psi: (0,1) \rightarrow (- \infty, \infty)$ is a bijective mapping.  We shall assume that
\be{
\label{map_endpts}
\psi(0) = -\infty,\qquad \psi(1) = \infty,
}
and also that
\be{
\label{map_symmetry}
\psi^{-1}(s) + \psi^{-1}(-s) = 1,\quad s \in (-\infty,\infty).
}
The latter assumption is not strictly necessary.  But it is useful to simplify some of the arguments and will in practice be satisfied by all mappings considered in this paper.

\subsection{Variable transform methods with domain truncation}
Given $f(x)$, we use $F(s)$ to denote the transplant of $f(x)$ to the new $s$-variable:
\ea{
F(s) = f( \psi^{-1}(s) ), \quad s \in (- \infty, \infty).
}
Approximation of $F(s)$ is achieved via domain truncation.  Let $L > 0$ and define the new function on the interval $[-1,1]$ by
\ea{
\label{E:FL}
F_L(y) = F( L y ), \quad y \in [-1,1].
}
For $n \in \bbN$, let $P_{n,L}(y)$ be the approximation of $F_L(y)$.  As mentioned, Chebyshev interpolation was used in \cite{AdcockRichardsonMappings} for this task.  In this paper we shall instead use a Fourier-type approximation, which will be introduced in the next section.
With this in hand,  the final approximation to $f(x)$ over the interval $[0,1]$ is defined as follows:
\be{
\label{E:AproxPsi}  p_{n,L} (x) =  \left\{  \begin{array}{cl}  F_L(-1) & x \in [0,x_L) \\
P_{n,L}( \psi(x)/L ) & x \in [x_L, 1-x_L] \\
F_L(1) & x \in (1- x_L,1] \end{array} \right. .
}
Here we use the notation $x_L = \psi^{-1}(-L) = 1-\psi^{-1}(L)$.  For the parametrized mappings, indexed by a parameter $\alpha$, we will also write $p_{n,L} = p_{n,L,\alpha}$ to make this dependence explicit.  With this in hand, we note that the approximation error $\| f - p_{n,L} \|_{[0,1]} = \max_{x \in [0,1]} | f(x) - p_{n,L}(x) |$ is given by
\be{
\label{E:errAprox}
\| f - p_{n,L} \|_{[0,1]} = \max \{ \| F_L - P_{n,L} \|_{[-1,1]}, \|f - f(x_L) \|_{ [0,x_L] },  \|f - f(1-x_L) \|_{[1-x_L,1] } \}.
}
Throughout the paper we will refer to the first term as the \textit{interior} error and the latter two terms as the \textit{endpoint} errors.

\subsection{The exponential and double-exponential maps}
The standard exponential mapping $\psi_{E}$ and its inverse are defined by
\be{
\label{E:PsiE}
\psi_E(x) = \log\left(  \frac{x}{1-x} \right), \quad  \psi_E^{-1}(s) =  \frac{1 }{1+\exp(-s)}.
}
Similarly, the standard double-exponential mapping $\psi_{DE}$ and its inverse are given by
\be{
\label{E:PsiDE}
\psi_{DE}(x) = \asinh \left( \frac{1}{\pi} \log \left(   \frac{ x } { 1-x}  \right) \right)     , \quad   \psi^{-1}_{DE}(s) =   \frac{1} { 1+\exp( -\pi \sinh(s) ) } .
}
We refer to \cite{LundBowers,Richard13Thesis,StengerBook93,StengerBook10} for background on these mappings.

\subsection{The parametrized exponential map}
As discussed in \cite{AdcockRichardsonMappings}, the exponential map $\psi_E$ may be undesirable in practice, since it requires more analyticity of the function $f$ in the interior than at the endpoints $x=0$ and $x=1$.  In particular, this leads to the poor resolution properties for oscillatory functions mentioned previously.  This observation can be understood by looking at the image of the strip
\be{
\label{strip}
S_{\beta} = \left \{ z \in \bbC : | \Im(z) | < \beta \right \},
}
of half-width $\beta$ under $\psi^{-1}_{E}$.  The corresponding region $\psi^{-1}_E(S_{\beta})$ turns out to be lens-shaped and formed by two circular arcs meeting with half-angle $\beta$ at $x=0$ and $x=1$.

Introduced in \cite{AdcockRichardsonMappings}, the parametrized exponential map $\psi_{SE}$ seeks to overcome this issue by enforcing that the strip $S_{\alpha}$ be mapped to a more regularly-shaped region.  In particular, $\psi^{-1}_{SE}(S_{\alpha}) = \tilde{S}_{\alpha}$ is a so-called \textit{two-slit strip} region of half-width $\alpha$:
\be{
\label{slitstrip}
\tilde{S}_{\alpha} = S_{\alpha} \backslash \{ (-\infty , 0 ] \cup [1,\infty) \}.
}
Such a transformation can be constructed via the Schwarz--Christoffel formula \cite{DriscollTrefethenSchwarzChristoffel}, leading to the following forms for $\psi_{SE}$ and its inverse:
\be{
\begin{aligned}
\label{E:PsiSE}
\psi_{SE}(x; \alpha)  &=  \frac{ \alpha} { \pi } \log \left( \frac{   \exp( \pi x/\alpha ) - 1 }  { 1 -  \exp( \pi (x-1) /\alpha )   } \right)
\\
 \psi^{-1}_{SE}(s; \alpha) &=  \frac{ \alpha} { \pi } \log \left( \frac{  1+  \exp( \pi (s+1/2)/\alpha ) }  { 1 +  \exp( \pi (s-1/2) /\alpha ) } \right ).
\end{aligned}
}
Note that the parameter $\alpha > 0$ is user-determined, and is chosen in conjunction with $L>0$ so as to obtain the best accuracy and resolution power.     
Later we will show that optimal choices for these parameters -- in the sense that they lead to the smallest error bound -- are given by
\bes{
\alpha = \alpha_0 / \sqrt{n},\qquad L = 1/2 + L_0,
}
for constants $\alpha_0,L_0 > 0$.  We refer to \S \ref{s:psiSE} for details.

\subsection{The parametrized double-exponential map}
The new map we introduce in this paper is a parametrized double-exponential map $\psi_{SDE}$.  It is defined via its inverse as follows:
\be{
\label{E:PsiSDE}
 \psi^{-1}_{SDE}(s; \alpha)  =
 \frac{ \alpha} { \pi } \log \left( \frac{  1+  \exp \left( \pi \left(   s  +1/2 \right)/\alpha + \frac{ \sinh(\pi s/ \alpha) }{\cosh( \pi /(2 \alpha) )}   \right) }  { 1 +  \exp \left( \pi \left(    s  -1/2 \right) /\alpha + \frac{ \sinh(\pi s/ \alpha) }  { \cosh( \pi /(2 \alpha) )}  \right) } \right).
}
Much as the standard double-exponential map $\psi_{DE} = g \circ \psi_E$ is a composition of the exponential map $\psi_E$ with the function $g(t) = \asinh(t/\pi)$, the parametrized double-exponential map $\psi_{SDE} = g \circ \psi_{SE}$ is a composition of the parametrized exponential map $\psi_{SE}$ and the function $g = g(t ; \alpha)$, defined through its inverse by
\be{
\label{E:gfunction}
g^{-1} (t; \alpha) = t + \frac{ \alpha \sinh(\pi t/ \alpha) }  { \pi \cosh( \pi /(2 \alpha) )}.
}
Note that $\psi_{SDE}$ involves a parameter, $\alpha >0$, which, as with $\psi_{SE}$, will be chosen to ensure best accuracy and resolution power.  Specifically, we make the following choice of parameters:
 \bes{
\alpha = L_0 \pi / (\pi/2 + W(c n)),\qquad  L = L_0 + 1/2.
}
Here $L_0, c>0$ are constants and $W$ is the Lambert-W function.  See \S \ref{s:psiSDE} for details.

Unfortunately $\psi_{SDE}(x ; \alpha)$, or more specifically, $g(t;\alpha)$, does not have an explicit form.  Whilst this does not impact computation of the approximation $p_{n,L,\alpha}$, it does affect one's ability to evaluate $p_{n,L,\alpha}$ at an arbitrary point $x$.  However, practical computation can be achieved via a few steps of Newton's method, for example.

Much like for $\psi_{SE}$, the specific choice for \R{E:gfunction} is motivated by the desire to preserve a certain amount of ``strip-behaviour'' near the real line.  The practical derivation of \R{E:gfunction} is also achieved using the Schwarz--Christoffel approach.  We refer to Appendix \ref{supp1} for the details.

\subsection{Resolution power}
In this paper, we not only derive error bounds and convergence rates for the various methods, we also investigate their ability to resolve oscillatory behaviour.  This is a traditional topic in the study of numerical algorithms \cite{BADHFEResolution,boyd,naspec}, and is usually done by studying the  complex exponential $f(x) = \E^{2 \pi \I \omega x}$.  This strategy has the benefit of providing a very clear quantitative measure of a numerical scheme -- the number of \textit{points-per-wavelength} (ppw) required to resolve an oscillatory function -- and therefore provides a direct way of comparing different methods.

Let $\{ \Psi_n \}_{n \in \bbN}$ be a sequence of numerical approximation with $\Psi_n$ having $n$ degrees of freedom.  For $\omega \geq 0$ and $0 < \delta < 1$ we define the $\delta$-resolution of $\{ \Psi_n \}_{n \in \bbN}$ as
\be{
\label{E:deltaresolution}
\cR(\omega ; \delta) = \min \left \{ n \in \bbN : \| \E^{-2 \pi \I \omega \cdot} - \Psi_n(\E^{-2 \pi \I \omega \cdot} ) \|_{[0,1]} \leq \delta \right \}.
}
We say the method $\{ \Psi_n \}_{n \in \bbN}$ has \textit{linear} resolution power if $\cR(\omega ; \delta) = \ord{\omega}$ as $\omega \rightarrow \infty$ for any fixed $0 < \delta < 1$ and \textit{sublinear} resolution power otherwise.  In the case of the former, we define the \textit{resolution constant} (the points-per-wavelength value) as
\be{
\label{E:resolutionconstant}
r = \limsup_{\omega \rightarrow \infty} \{ \cR(\omega ; \delta ) / \omega \}.
}
Note that $r$ need not be independent of $\delta$, but this will be the case for all methods considered in this paper.

For classical Chebyshev interpolation the resolution constant is $r = \pi$.  Conversely, Fourier interpolation has $r = 2$, provided the $\omega$ in \R{E:resolutionconstant} are restricted to integer values.  Thus, although both schemes have linear resolution power, Fourier interpolation is more efficient by a factor of $\pi/2$ for resolving periodic oscillations. Such an improvement is what we shall seek to achieve in this paper when using variable transform methods for approximating functions with singularities.\footnote{Although the numerical methods in this paper are designed for functions with endpoint singularities, we formulate \R{E:deltaresolution} in terms of the entire functions $f(x) = \E^{-2 \pi \I \omega x}$.  However, our results are unchanged if we allowed the more general form $f(x) = g(x) \E^{-2 \pi \I \omega x}$, where $g(x)$ is analytic on $(0,1)$, continuous on $[0,1]$, and independent of $\omega$.  For simplicity, we merely consider the case $g(x) = 1$.
}

\section{Approximation strategy in the truncated domain}\label{s:FourierLike}
As discussed, the function $F(s)$ converges rapidly to its limiting values as $|s| \rightarrow \infty$.  For large $L$, the normalized function $F_L(y)$ is therefore flat near the endpoints $y=\pm 1$.  One option to approximate $F_L(y)$ efficiently would be to first subtract its endpoint values with a linear function of $y$ and then expand the remainder in a Fourier series.  However, while theoretically sound, this may cause practical issues, especially when incorporating variable transform techniques into numerical computing packages (e.g.\ \texttt{Chebfun}) \cite{Richard13Thesis}.

\subsection{Cosine expansions}
Instead, we shall approximate $F_L(y)$ using a cosine expansion.  The advantage of this expansion is that it retains the key properties of Fourier expansions vis-a-vis resolution power -- that is, a factor of $\pi/2$ better than Chebyshev expansion -- without the requirement of periodicity.  In particular, such an expansion (also known as a modified Fourier expansion \cite{BAthesis,MF1}) is uniformly convergent for continuously-differentiable functions without a periodicity assumption \cite{BA1}.

This expansion can be written as
\bes{
F_L(y) = \sum^{\infty}_{k=0} c_k \cos (k \pi (y+1)/2),\qquad y \in [-1,1],
}
where
\be{  \label{E:FC}
 c_0 = \frac12 \int^{1}_{-1} F_L(y) \D y,\qquad c_k = \int^{1}_{-1} F_L(y) \cos (k \pi (y+1)/2) \D y,\quad k =1,2,\ldots,
}
are the coefficients of $F_L$.  For a given $n \in \bbN$, we define the discrete coefficients by
 \bes{
\tilde{c}_k =  \frac{  2 \gamma_{k} }{n} \sum^{n}_{j=0} \gamma_{j} F_L(y_j) \cos(j k \pi / n),\qquad k=0,\ldots,n,
}
where $y_{j} = -1+2j/n$, $j=0,\ldots,n$ and $\gamma_0 = \gamma_n = 1/2$ and $\gamma_{k} = 1$ otherwise.  Note that these coefficients can be computed in $\ord{n \log n}$ operations using FFTs.  We now obtain the following approximation to $F_L$:
\bes{
P_{n,L}(y) = \sum^{n}_{k=0} \tilde{c}_k \cos(k \pi(y+1)/2).
}
Observe that if $F_L$ is Lipschitz continuous on $[-1,1]$ then we have the following absolutely convergent expression for the error
\be{
\label{error_aliasing}
F_L(y) - P_{n,L}(y) = \sum_{k > n} c_k \left ( \cos (k \pi (y+1)/2) - \cos (k' \pi (y+1)/2) \right ),
}
where $k' = | \mathrm{mod}(k+n-1,2n) - (n-1) |$ \cite{TrefethenATAP}.

\begin{remark}
\label{r:whynotsinc}
Domain truncation followed by approximation is not the only strategy for variable transform methods.  A common alternative would be $\mathrm{sinc}$ expansion on the real line.  However, it has been argued in \cite{RichardsonSincFun} that this may not be best solution in practice, especially in numerical computing environments where the domain truncation $L$ is typically fixed before approximation.  As we show in this paper, the numerical method based on cosine expansions achieves exactly the same convergence rate and resolution power as the corresponding $\mathrm{sinc}$-based method.
\end{remark}

\subsection{Error analysis}
The downside of the cosine-based method over a Chebyshev-based method is that the analysis is more involved, since the function $F_L(y)$ is not exactly flat at the endpoints $y = \pm 1$.  To this end, we now present the following result, which will be crucial later in deriving error bounds for all maps under consideration:

\lem{
\label{lemConvPsiGeneral}
Let $\psi : (0,1) \rightarrow (-\infty,\infty)$ be an invertible mapping satisfying \R{map_endpts} and \R{map_symmetry}, and let $f$ be analytic and bounded in $\psi^{-1}(S_{\beta})$ for some $\beta > 0$, where $S_\beta$ is as in \R{strip}.  Define
\be{
\label{E:M1}
M_\psi = M_\psi(\beta ; f) = \sup_{z \in S_\beta} | f(\psi^{-1}(z)) | < \infty,
}
and suppose that there exists $0 < \tau < 1$ such that
\ea{
\label{E:M2}
N_\psi = N_\psi (\beta,L ; f) =   \sup_{|z\mp L| \leq \beta} |f(\psi^{-1}(z)) -f(\iota_{\pm}) | / |\iota_{\pm}-\psi^{-1}(z)|^{\tau}  < \infty,
}
where $\iota_{+} = 1$ and $\iota_{-} = 0$,
and let $p_{n,L}$ be given by \R{E:AproxPsi}.  If $L \leq n$ then
\be{  \label{E:lemma}
\| f - p_{n,L} \|_{[0,1]} \leq  3 N_{\psi} \bar{\beta}^{-1} (C_\psi)^{\tau}+114  M_{\psi} \bar{\beta}^{-2} n \E^{-\beta n \pi/(2L)} ,
}
for any $n \in \bbN$, where $\bar{\beta} = \min \{ \beta, 1\}$ and
\bes{
C_{\psi} = C_{\psi}(\beta,L) =  \sup_{|s+L| \leq \beta} | \psi^{-1}(s) |.
}
}
Note that the second term in \R{E:lemma}, proportional to $M_{\psi} \E^{-\beta n \pi / (2L)}$, corresponds to the interior error $\| F_L - P_n \|_{[-1,1]}$ in \R{E:errAprox}, whereas the first, proportional to $N_{\psi} (C_{\psi})^{\tau}$, corresponds to the endpoint errors $\| f - f(x_L) \|_{[0,x_L]}$ and $\| f - f(1-x_{L}) \|_{[1-x_L,1]}$.
We remark also that the condition $L \leq n$ is not fundamental, and could be relaxed at the expense of a more complicated statement.  We include it since it leads to a simpler error bound and is always satisfied in practice for the mappings we use in this paper.

\begin{proof}[Proof of Lemma \ref{lemConvPsiGeneral}]
We shall bound the three terms in \R{E:errAprox} separately.    Consider the first term.  By \R{error_aliasing}, we have
\be{
\label{GLPn_est}
\| F_L - P_{n,L} \|_{[-1,1]} \leq 2 \sum_{k>n} | c_k |.
}
Hence, we now wish to estimate the coefficients $c_k$.  Let $t \geq 1$.    After integrating the formula \R{E:FC} by parts $2t$ times we find that
\eas{
c_k =& \sum^{t-1}_{r=0} \frac{(-1)^{r}}{(k \pi/2)^{2r+2}} \left ( (-1)^k F^{(2r+1)}_L(1) - F^{(2r+1)}(-1) \right )
\\
&+ \frac{(-1)^{t}}{(k \pi/2)^{2t}} \int^{1}_{-1} F^{(2t)}_{L}(y) \cos(k \pi(y+1)/2) \D y.
}
Thus
\be{
\label{ck_bound1}
|c_k| \leq 2 \sum^{t-1}_{r=0} \frac{\max \{ | F^{(2r+1)}_{L}(1) | , | F^{(2r+1)}_L(-1) | \}}{(k\pi/2)^{2r+2} } + \frac{2}{(k \pi/2)^{2t}} \| F^{(2t)}_{L} \|_{[-1,1]}.
}
We now wish to estimate $|F^{(s)}_{L}(a)|$ for $-1 \leq a \leq 1$ and $s=0,1,2,\ldots$.  By Cauchy's integral formula
\bes{
F^{(s)}_{L}(a) = \frac{s!}{2 \pi \I} \oint_{C} \frac{ F_L(y) } { (y-a)^{s+1} } dy,
}
where $C$ is any circular contour of radius $0 < \rho < \beta / L$ centred at the point $y=a$.  Note that this contour lies within the strip $S_{\beta/L}$ and therefore within the region of analyticity of $F_L$.  By a change of variables, we find that
\bes{
F^{(s)}_{L}(a) = \frac{s! L^{s}}{2 \pi \I} \oint_{C} \frac{F(z)}{(z-L a)^{s+1}} \D z,
}
where now $C$ is any circular contour of radius $0 < \rho < \beta$ centred at the point $z = L a$.  Therefore
\be{
\label{F_L_deriv}
| F^{(s)}_{L}(a)| \leq \frac{s! L^{s}}{\rho^{s}} \sup_{ | z- L a | = \rho} | F(z) |.
}
Recall that $F(z) = f(\psi^{-1}(z))$.  In particular, we deduce that
\bes{
\| F^{(2t)}_{L} \|_{[-1,1]} \leq \frac{(2t)! L^{2t}}{\beta^{2t}} M_{\psi}.
}
Now suppose that $a = \pm 1$ in \R{F_L_deriv}.  Then,
\eas{
| F^{(2r+1)}_{L}(\pm 1) | &\leq \frac{(2r+1)! L^{2r+1}}{\beta^{2r+1}} \sup_{|z \mp L | = \beta} | f(\psi^{-1}(z))|
\\
 & \leq \frac{(2r+1)! L^{2r+1}}{\beta^{2r+1}} N_{\psi} \sup_{|z\mp L|=\beta} | \psi^{-1}(z) - \iota_{\pm} |^{\tau}.
}
By \R{map_symmetry}, we notice that
\bes{
\sup_{|z-L|=\beta} | \psi^{-1}(z) - 1 | = \sup_{|z+L| = \beta} | \psi^{-1}(z) | = C_{\psi}.
}
Substituting this and \R{F_L_deriv} into \R{ck_bound1} gives
\be{
\label{c_k_bound}
|c_k| \leq 2 N_{\psi} (C_{\psi})^{\tau} \sum^{t-1}_{r=0} \frac{(2r+1)! L^{2r+1}}{ (k \pi /2)^{2r+2} \beta^{2r+1} } + 2 M_{\psi} \frac{(2t)! L^{2t} }{(k \pi/2)^{2t}\beta^{2t}}.
}
and we now substitute this into \R{GLPn_est} to get
\eas{
\| F_L - P_{n,L} \|_{[-1,1]} \leq & 4 N_{\psi} (C_{\psi})^{\tau} \sum^{t-1}_{r=0}  \frac{(2r+1)! L^{2r+1}}{(\pi/2)^{2r+2} \beta^{2r+1}} \sum_{k > n} \frac{1}{k^{2r+2}}
\\
& + 4 M_{\psi} (2t)! \left ( \frac{2 L}{\beta \pi} \right )^{2t} \sum_{k>n} \frac{1}{k^{2t}}.
}
Using the fact that $\sum_{k > n} k^{-s-1} \leq s^{-1} n^{-s}$, we deduce the following:
\be{
\label{alamo}
\| F_L - P_{n,L} \|_{[-1,1]} \leq \frac{8 N_{\psi} (C_{\psi})^{\tau}}{\pi } \sum^{t-1}_{r=0}  (2r)! \left ( \frac{2L}{\beta \pi n} \right )^{2r+1}+ 4 M_{\psi} \frac{(2t)!}{2t-1} \left ( \frac{2L}{\beta \pi} \right )^{2t} n^{1-2t}.
}
We are now in a position to choose $t$.  For this, we consider two cases:

\vspace{1pc}
\noindent Case (i): $\beta n \pi / L \geq 4$.  Let $t = \lfloor \beta n \pi / (4 L) \rfloor \geq 1$.  Then
\bes{
\sum^{t-1}_{r=0}  (2r)! \left ( \frac{2L}{\beta \pi n} \right )^{2r} \leq \sum^{t-1}_{r=0} (2r)! (2t)^{-2r}.
}
We estimate this using the upper bound in Stirling's formula $n! \leq \sqrt{ 2 \pi } n^{n+1/2} \E^{-n} \E^{1/12}$.  This gives
\eas{
\sum^{t-1}_{r=0}  (2r)! \left ( \frac{2L}{\beta \pi n} \right )^{2r} & \leq 1 +\sum^{t-1}_{r=1} \sqrt{2\pi} \E^{1/12} (2r)^{2r+1/2} \E^{-2r} (2t)^{-2r}
\\
& \leq 1 + \sqrt{2\pi} \E^{1/12} \sum^{t-1}_{r=1}\sqrt{2r} \E^{-2r}.
}
Since the function $\sqrt{2x} \E^{-2x} $ is decreasing whenever $x \geq 1$ we have
\bes{
\sum^{t-1}_{r=1}  \sqrt{2r} \E^{-2r} \leq  \sum^{ \infty }_{r=1}  \sqrt{2r} \E^{-2r} \leq \int^{ \infty }_{0} \sqrt{2x} \E^{-x} \D x = \sqrt{\pi}/4,
}
and therefore
\bes{
\sum^{t-1}_{r=0}  (2r)! \left ( \frac{2L}{\beta \pi n} \right )^{2r} \leq 1 + \sqrt{2} \pi \E^{1/12} / 4.
}
This now gives
\bes{
\| F_L - P_{n,L} \|_{[-1,1]} \leq 3 N_{\psi} (C_\psi)^{\tau} +4 n M_{\psi} \frac{(2t)!}{2t-1} \left ( \frac{2L}{\beta  n \pi} \right )^{2t}.
}
We can use Stirling's formula once more to estimate the second term:
\eas{
\frac{(2t)!}{2t-1} \left ( \frac{2L}{\beta n \pi} \right )^{2t} \leq & \frac{1}{2t-1} \sqrt{2\pi} (2t)^{2t+1/2} \E^{-2t} \E^{1/12} \left ( \frac{2L}{\beta n \pi} \right )^{2t}
\\
=& \sqrt{2 \pi} \E^{1/12} \frac{\sqrt{2t}}{2t-1} \E^{-2t} \left ( \frac{4 L t}{\beta n \pi} \right )^{2t}
\\
\leq & 2 \sqrt{\pi} \E^{1/12} \E^{-2t} \leq 2 \sqrt{\pi} \E^{25/12} \E^{-\beta n \pi / (2L)}
}
Substituting this into the earlier expression now gives
\bes{
\| F_L - P_{n,L} \|_{[-1,1]} \leq 3 N_{\psi} (C_\psi)^{\tau} +8 \sqrt{\pi} \E^{25/12}  M_{\psi} n \E^{-\beta n \pi/(2L)}.
}
To complete the proof for this case, it therefore remains to estimate the other two terms in \R{E:errAprox}.  For $x \in [0,x_L]$, notice that
\bes{
|f(x) - f(x_L)| \leq | f(x) - f(0) | + | f(0) - f(x_L) | \leq 2 N_{\psi} | x_L |^{\tau}.
}
However, $|x_L| = |\psi^{-1}(-L)| \leq C_{\psi}$.  Therefore $\| f - f(x_L) \|_{[0,x_L]} \leq 2 N_{\psi} (C_{\psi})^{\tau}$ and we deduce that
\bes{
\| f - p_{n,L} \|_{[0,1]} \leq 3 N_{\psi} (C_\psi)^{\tau} +8 \sqrt{\pi} \E^{25/12}  M_{\psi} n \E^{-\beta n \pi/(2L)},\qquad \beta n \pi / L \geq 4.
}
Case (ii): $0 < \beta n \pi / L < 4$. First, set $t=1$ so that, by \R{alamo} and the fact that $L \leq n$, we have
\eas{
\| F_L - P_{n,L} \|_{[-1,1]} &\leq \frac{16}{\pi^2} N_{\psi} (C_{\psi})^{\tau} \frac{L}{\beta n} + \frac{32}{\pi^2} M_{\psi} \frac{L^2}{\beta^2 n}
\\
& \leq \frac{16}{\pi^2} N_{\psi} (C_{\psi})^{\tau} \bar{\beta}^{-1} + \frac{32 \E^2}{\pi^2} M_{\psi} \bar{\beta}^{-2} n \E^{-\beta n \pi / (2 L) }.
}
Since $\bar{\beta} \leq 1$ and the endpoint error is the same as in case (i), we get
\bes{
\| f - p_{n,L} \|_{[0,1]} \leq 2 N_{\psi}  \bar{\beta}^{-1} (C_{\psi})^{\tau} +  \frac{32 \E^2}{\pi^2} M_{\psi} \bar{\beta}^{-2} n \E^{-\beta n \pi / (2 L)},\qquad 0 < \beta n \pi / L < 4.
}
To complete the proof, we note that $32 \E^2 / \pi^2 \leq 8 \sqrt{\pi} \E^{25/12} < 114$ and then combine the estimates from cases (i) and (ii).
\end{proof}

\section{The exponential map $\psi_E$  }\label{s:psiE}
In this and the next three sections we analyze the convergence and resolution power for the numerical methods based on cosine expansions for the four maps $\psi_E$, $\psi_{SE}$, $\psi_{DE}$ and $\psi_{SDE}$.  In doing so, we also determine appropriate choices for the truncation parameter $L$, as well as the mapping parameter $\alpha$ (where appropriate).  We commence with $\psi_E$.

\subsection{Convergence}
Our main result for $\psi_E$ is as follows:
\thm{
\label{t:psiEconv}
Let $\psi = \psi_E$ be the mapping given by \R{E:PsiE}.  Suppose that $f$ is analytic and bounded in the domain $\psi^{-1}(S_{\beta})$ for some $0<\beta<\pi$.  Let $M_{\psi}$ be as in \R{E:M1} and suppose that there exists $0 < \tau < 1$ such that \R{E:M2} holds with constant $N_{\psi}$.  Let $p_{n,L}$ be the approximation defined by \R{E:AproxPsi}.  If $L = c \sqrt{n}$ for some $c>0$ then
\be{
\label{psiEerr}
\| f - p_{n,L} \|_{[0,1]} \leq A \left [  M_{\psi} \beta^{-2} n \exp(-\beta \pi \sqrt{n}/(2c)) + N_{\psi} \beta^{-1} \exp(-\tau c \sqrt{n}) \right ],
}
for all $n \geq c^{-2}(\pi +\log(2))^2$, where the constant $A \leq 114 \pi^2 \approx 1125$.}
\prf{
We shall use Lemma \ref{lemConvPsiGeneral}.  First, we recall from \cite{AdcockRichardsonMappings} that for $0 < \beta < \pi$, $\psi^{-1}(S_{\beta})$ is well-defined and corresponds to a lens-shaped region formed by circular arcs meeting with half-angle at $x=0$ and $x=1$.  Also, for $0 < \beta < \pi$ we have $\bar{\beta} = \min \{ \beta , 1 \} \geq \beta / \pi$.  Hence
\bes{
\| f - p_{n,L} \|_{[0,1]} \leq 3 \pi N_{\psi} \beta^{-1} (C_{\psi})^{\tau} + 114 \pi^2 M_{\psi} \beta^{-2} n \E^{-\beta \pi \sqrt{n} / (2 c) }.
}
It remains to estimate the constant $C_{\psi}$.  We have
\bes{
C_{\psi} =  \sup_{|s+L|\leq \beta} | \psi^{-1}_{E}(s)  |  =  \sup_{|s+L|\leq \beta} \left|  \frac{\E^{s}}{\E^{s}+1 }  \right| =  \sup_{|s+L|\leq \beta} \frac{1}{| 1 + \E^{-s} |} \leq \left ( \E^{L-\beta}-1 \right )^{-1}.
}
Thus $C_{\psi} \leq 2 \E^{\beta - c \sqrt{n}}$ for all $n \geq c^{-2} (\pi + \log(2))^2 \geq c^{-2} (\beta + \log(2))^2$, as required.  In particular,
\bes{
(C_{\psi})^{\tau} \leq (2 \E^{\beta})^{\tau} \E^{-\tau c \sqrt{n}} \leq 2 \E^{\pi} \E^{-\tau c \sqrt{n}}.
}
Substituting this in to the previous bound now gives the result.
}

Note that the choice $L = c \sqrt{n}$ is made here to ensure that the two terms in the error expression \R{E:lemma} decay at the same, root-exponential, rate.  Allowing $L$ to scale either faster or slower with $n$ would lead to a slower convergence rate.

\begin{figure}
\begin{center}
\begin{tabular}{c c}
\includegraphics[width=5.5cm]{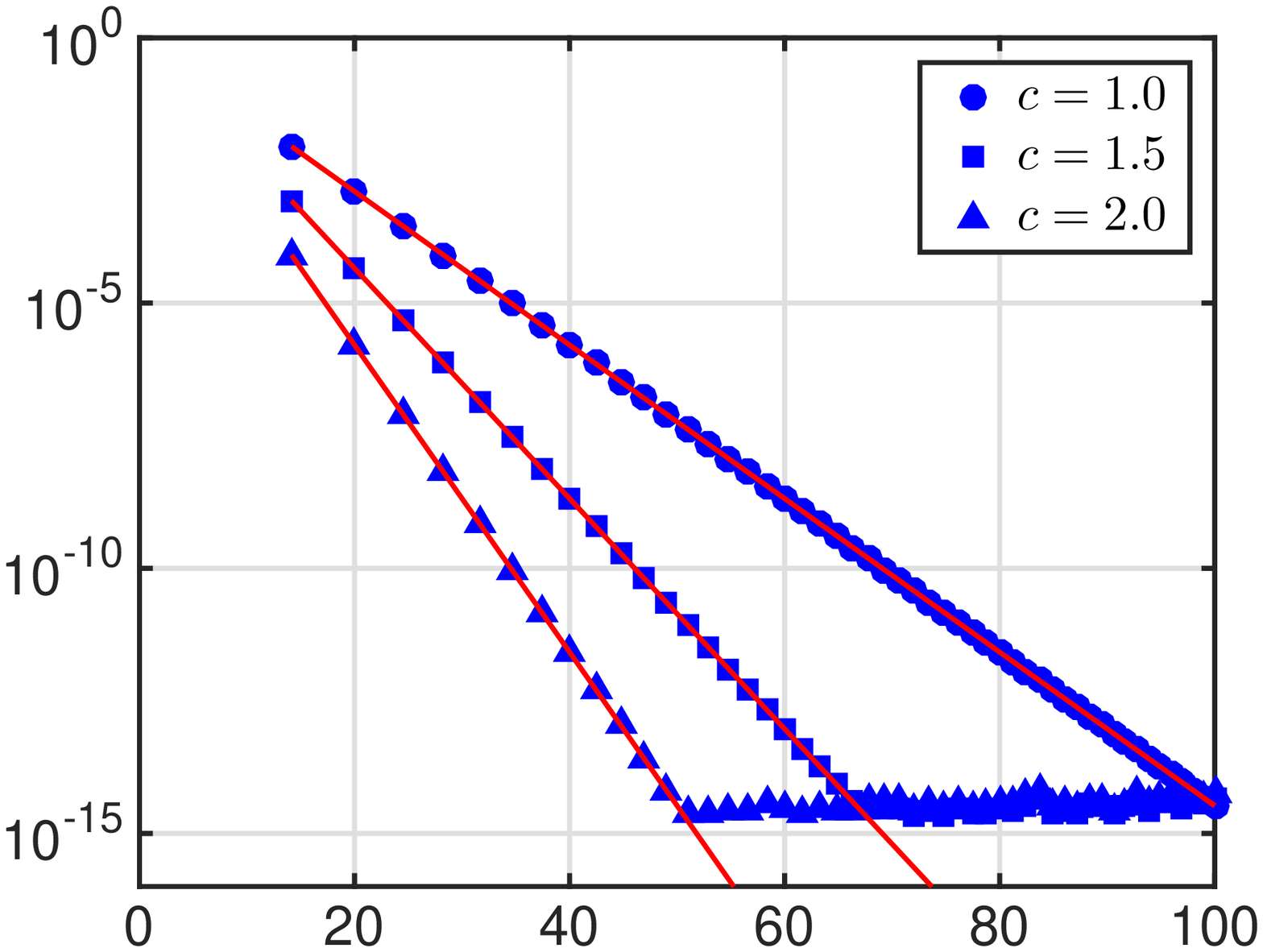} &\includegraphics[width=5.5cm]{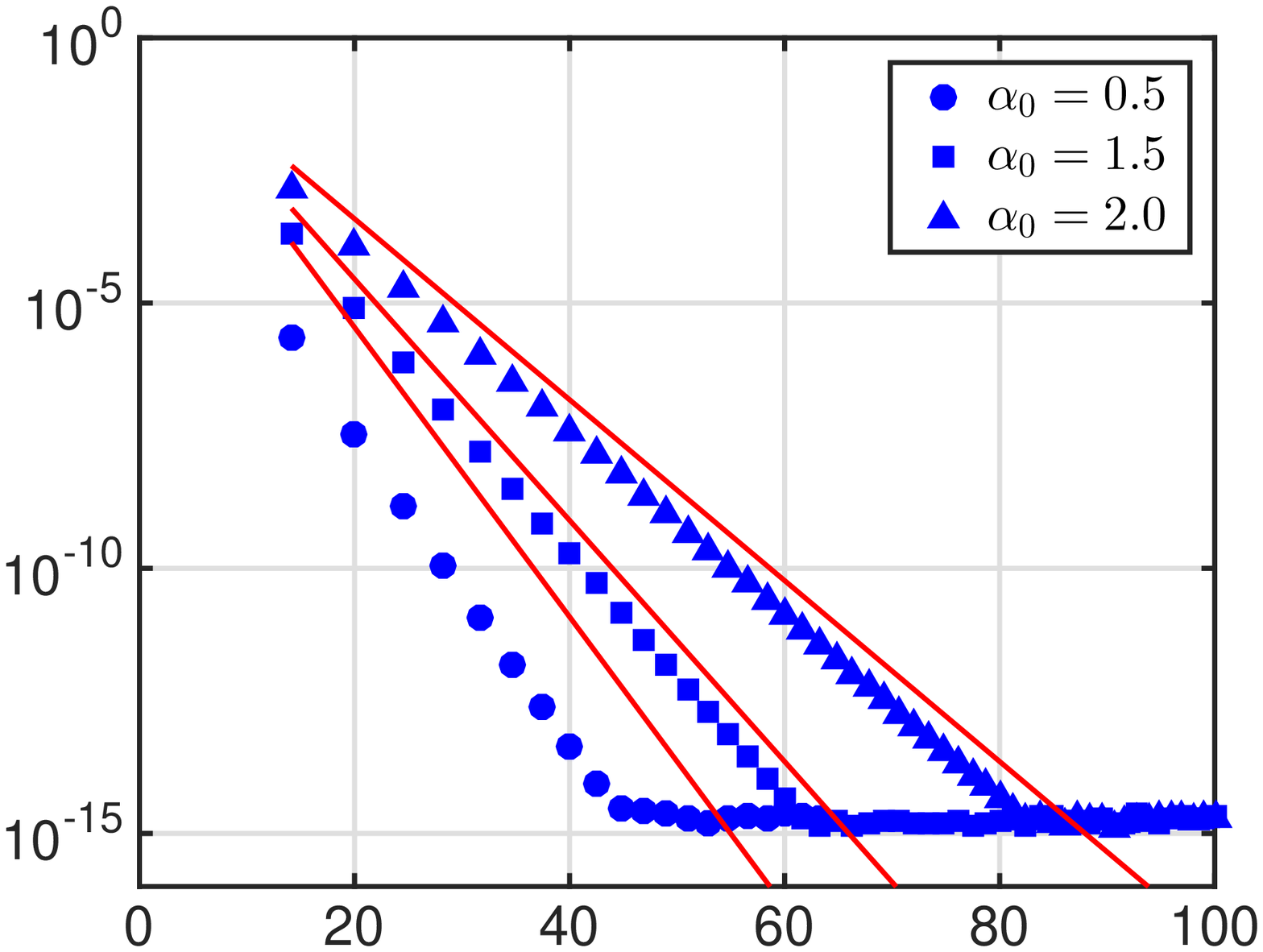}
\\
{\small $\psi_{E}$, $L = c \sqrt{n}$} & {\small $\psi_{SE}$, $L = 1.25$, $\alpha = \alpha_0 / \sqrt{n}$}
\end{tabular}
\end{center}
\caption{The error against $\sqrt{n}$ for $f(x) = x^{1/3}$.  The red lines shows the quantities $\rho^{-\sqrt{n}}$, where $\rho$ is given by \R{rho_psiE} or \R{rho_psiSE} for $\psi_E$ or $\psi_{SE}$ respectively with $\tau = 1/3$.}
\label{f:PsiEConv}
\end{figure}

In Fig.\  \ref{f:PsiEConv} we plot the error $\| f - p_{n,L} \|_{[0,1]}$ for several choices of $c$.\footnote{In this and all subsequent experiments in this paper we compute this error using \R{E:errAprox}.  The various max-norms are replaced by a discrete maximum over a grid of $20,000$ equally-spaced points.}
  As with all other experiments in this paper, the results shown in this figure were computed in \textit{Matlab} using double precision.  Theorem \ref{t:psiEconv} predicts that the error decays root-exponentially fast in $n$ with index
\be{
\label{rho_psiE}
\rho = \min \{ \exp(\beta \pi / (2c)) , \exp(\tau c ) \},
}
that is $\| f - p_{n,L} \|_{[0,1]} =\ordu{n \rho^{-\sqrt{n}}}$ as $n \rightarrow \infty$, and this is in good agreement with the results shown in this figure.

\rem{
\label{r:pi_over_two}
The contribution from the interior error term $\beta \pi / (2c)$ is identical to the corresponding estimate for $\mathrm{sinc}$ interpolation (see \cite[Thm 2.2]{Richard13Thesis}), and precisely $\pi/2$ times larger than the estimate for Chebyshev interpolation (see \cite[Thm.\ 3.1]{AdcockRichardsonMappings}).  This is a well-known phenomenon when comparing Fourier and Chebyshev interpolation for analytic functions \cite{HaleTrefethenMapped} and one which perhaps unsurprisingly carries over to the case of variable transform methods for approximating functions with singularities.  Note that the endpoint contribution $\exp(\tau c)$ is identical to that of the $\mathrm{sinc}$ and Chebyshev-based methods, as one would expect, since it depends only on the map $\psi$ and not on the approximation scheme used.
}

\subsection{Resolution power}
To analyze the resolution properties of  $\psi_{E}$ (and the other maps in the paper), we proceed by estimating the quantities $M_{\psi}$ and $N_{\psi}$ for the function $f(x) = \E^{-2 \pi \I \omega x}$.  This leads to the following result:

\thm{
\label{t:PsiE_respower}
Let $\psi_E$ be the mapping given by \R{E:PsiE} and suppose that $f(x) = \exp(-2 \pi \I \omega x)$ for some $\omega \geq \pi+\log(2)$.  If $p_{n,L}$ is the approximation defined by \R{E:AproxPsi} and $L = c \sqrt{n}$ for some $c>0$ then, for $n > n^* = c^2 \omega^2$,
\ea{
\label{resn_err}
\| f - p_{n,L} \|_{[0,1]} \leq A \Big [ & \left(1-(n^*/n)^{1/4} \right )^{-1} \exp \left ( -\pi \omega H\left(\sqrt{n/n^*}\right)\right )
\\
& + 2\pi \omega  \left(1-(n^*/n)^{1/4} \right )^{-2} \exp \left ( 4 \pi \omega \E^{\pi-c \sqrt{n}} - c \sqrt{n}\right ) \Big ], \nn
}
where the constant $A>0$ is as in Theorem \ref{t:psiEconv} and the function $H$ is given by $H(t) = 0$ for $0 \leq t < 1$ and $H(t) = t \arccos(1/\sqrt{t}) - \sqrt{t-1}$ for $t>1$.  In particular,
\be{
\label{E:PsiE_respower}
\limsup_{\omega \rightarrow \infty} \cR(\omega ; \delta ) / \omega^2 \leq c^2,\qquad 0 < \delta < 1.
}
}
\prf{
First, consider $M_{\psi}$.  By the maximum modulus principle, we have
\bes{
M_{\psi} = \sup_{\substack{z = x \pm \I \beta \\ x \in \bbR}} \left|\exp \left( - 2\pi  \I \omega \frac{ 1 }  { 1 + \E^{-z} } \right) \right | =  \exp \left (2 \pi \omega \sup_{\substack{z = x \pm \I \beta \\ x \in \bbR}} \Im \left ( \frac{1}{1+\E^{-z}} \right ) \right ).
}
Let $z = x + \I \beta$ and observe that
\bes{
\sup_{x \in \bbR} \Im \left(  \frac{ 1 }  { 1 + \E^{-x-i \beta} } \right) =  \sup_{x \in \bbR} \frac{ \sin(\beta) }  { \E^{x} + 2  \cos(\beta) + \E^{-x} } = \frac{\tan(\beta/2)}{2}.
}
After noting that the corresponding term with $z = x - \I \beta$ is always negative, we deduce that $M_{\psi} = \exp(\pi \omega \tan(\beta/2))$.  

Now consider $N_{\psi}$.  Letting $\tau = 1$ (since $f(x) = \exp(-2 \pi \I \omega x)$ is entire), we find that $N_{\psi} \leq 2 \pi \omega \sup_{|z \pm L| \leq \beta} | f(\psi^{-1}(z)) |$.  By similar arguments to those given above, we get
\eas{
N_{\psi}  \leq 2 \pi \omega \sup_{\substack{| x \pm L | \leq \beta \\ | y | \leq \beta} } \exp \left ( 2 \pi \omega \frac{\sin(y)}{\E^x + 2 \cos(y) + \E^{-x} } \right )
 \leq 2 \pi \omega \exp \left (2 \pi \omega / \left ( \E^{L-\beta} - 1 \right ) \right ).
}
The assumptions on $L$, $\beta$ and $n$ now give $N_{\psi} \leq 2 \pi \omega \exp \left ( 4 \pi \omega \E^{\pi-c \sqrt{n}}  \right )$.

Combining this with the estimate for $M_{\psi}$, we deduce from Theorem \ref{t:psiEconv} that $\| f - p_{n,L} \|_{[0,1]} $ is bounded by
\bes{
 A \left ( \beta^{-2} n \exp \left ( \pi \omega \tan(\beta/2) - \sqrt{n} \pi \beta / (2c) \right ) + 2\pi \omega \beta^{-1} \exp \left ( 4 \pi \omega \E^{\pi -c \sqrt{n}}- c \sqrt{n}\right ) \right ),
}
which holds for any $0 < \beta < \pi$.  Besides the factor $\beta^{-1}$, the second term is independent of $\beta$.  Hence, disregarding the $\beta^{-2}$ factor, we now optimize the first term with respect to $\beta$.  For $n \geq n^*$ the function
\bes{
\pi \omega \tan(\beta/2) - \sqrt{n} \pi \beta / (2c),\qquad 0 \leq \beta < \pi,
}
has a local minimum at $\beta = 2 \arccos ( (n^*/n)^{1/4})$ and takes the value $-\pi \omega H(\sqrt{n/n^*})$ there.  Substituting this into the previous bound and noticing that this choice of $\beta$ satisfies $\beta \geq \pi(1-(n^*/n)^{1/4})$ now gives \R{resn_err}.

For \R{E:PsiE_respower} we first let $n/n^* = 1+ k \omega^{-1/2} > 1$ and then consider the behaviour of \R{resn_err} as $\omega \rightarrow \infty$.  Note that the right-hand side of \R{resn_err} is $o(1)$ for this choice of $n$ as $\omega \rightarrow \infty$.  In particular, for each fixed $\delta$ we have $\| f - p_{n,L} \|_{[0,1]} \leq \delta$ for all large $\omega$.  This now gives the result.
}

\begin{figure}
\begin{center}
\begin{tabular}{ccc}
\includegraphics[width=5.5cm]{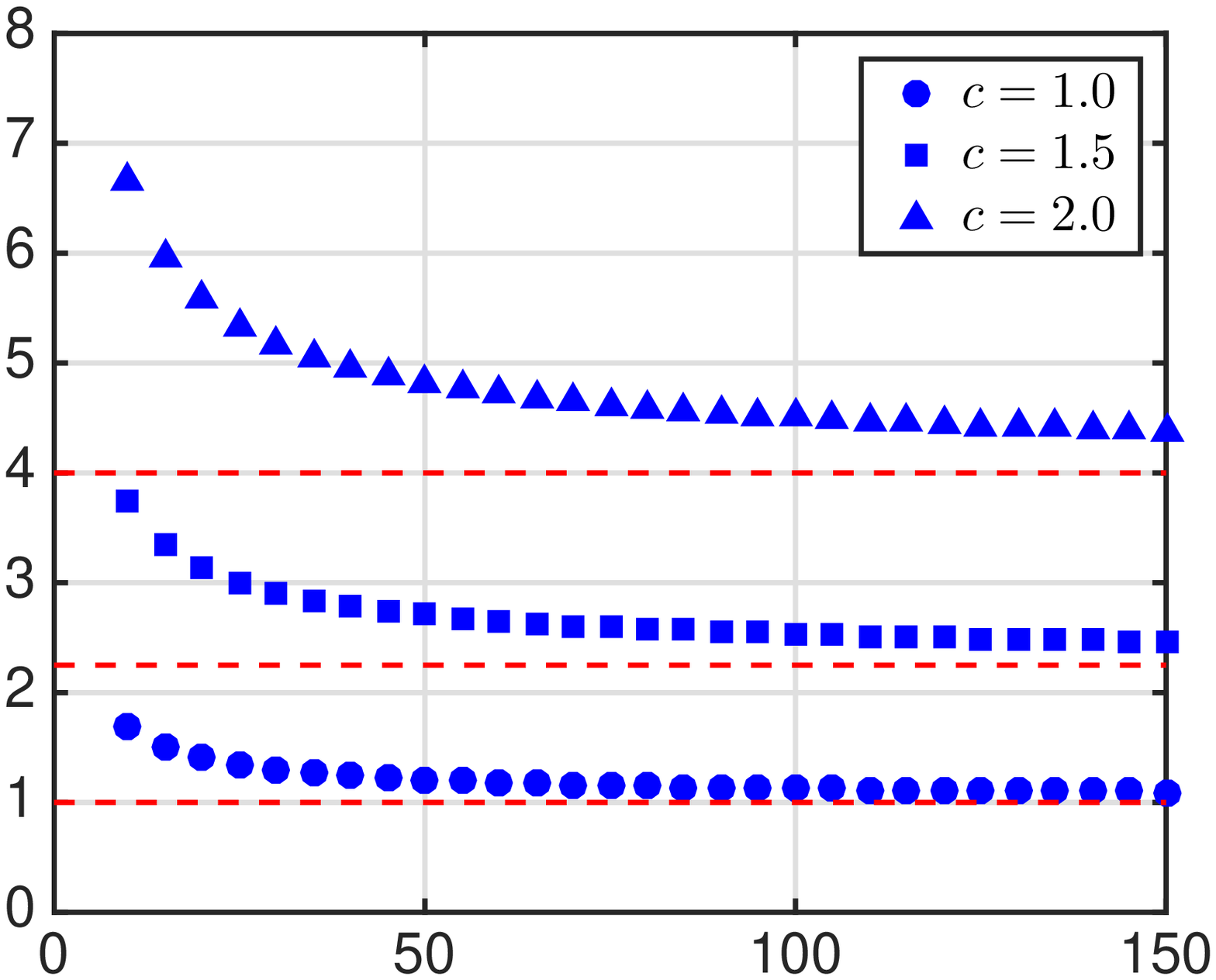} && \includegraphics[width=5.5cm]{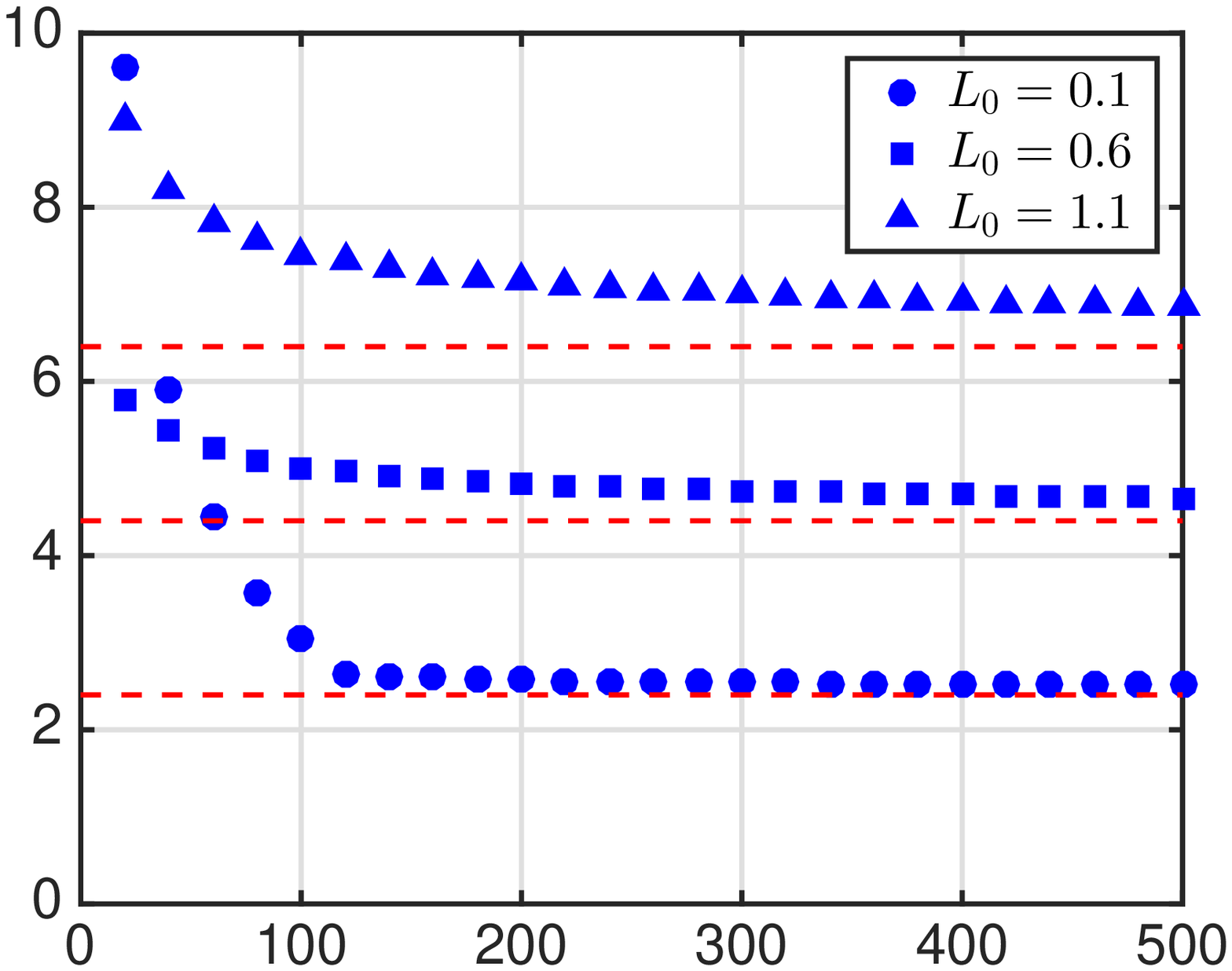}
\\
{\small $\psi_E$, $L = c \sqrt{n}$} && {\small $\psi_{SE}$, $L = L_0+1/2$, $\alpha = 0.8$}
\end{tabular}
\end{center}
\caption{Left: the quantity $\cR(\omega ; \delta) / \omega^2$ against $\omega$ with $\delta = 10^{-2}$ for $\psi_E$.  The dashed red line shows the theoretical resolution constant $c^2$. Right: the quantity $\cR(\omega ; \delta) / \omega$ with $\delta = 10^{-2}$ for $\psi_{SE}$.  The dashed red line shows the theoretical resolution constant $4(L_0+1/2)$.
}
\label{f:PsiEResolution}
\end{figure}

Similar to a result proved in \cite{AdcockRichardsonMappings} for Chebyshev approximation, this theorem shows that $\psi_E$ has sublinear resolution power, scaling quadratically with the frequency.  Numerical verification of \R{E:PsiE_respower} is shown in the left panel of Fig.\  \ref{f:PsiEResolution}.  Interestingly, the constant $c^2$ is precisely $(\pi/2)^2$ times smaller than the corresponding constant for Chebyshev approximation \cite[Thm.\ 3.2]{AdcockRichardsonMappings}, as one might expect given the differences in the approximation schemes (recall Remark \ref{r:pi_over_two}).

\begin{figure}
\begin{center}
\begin{tabular}{ccc}
\includegraphics[width=4.0cm]{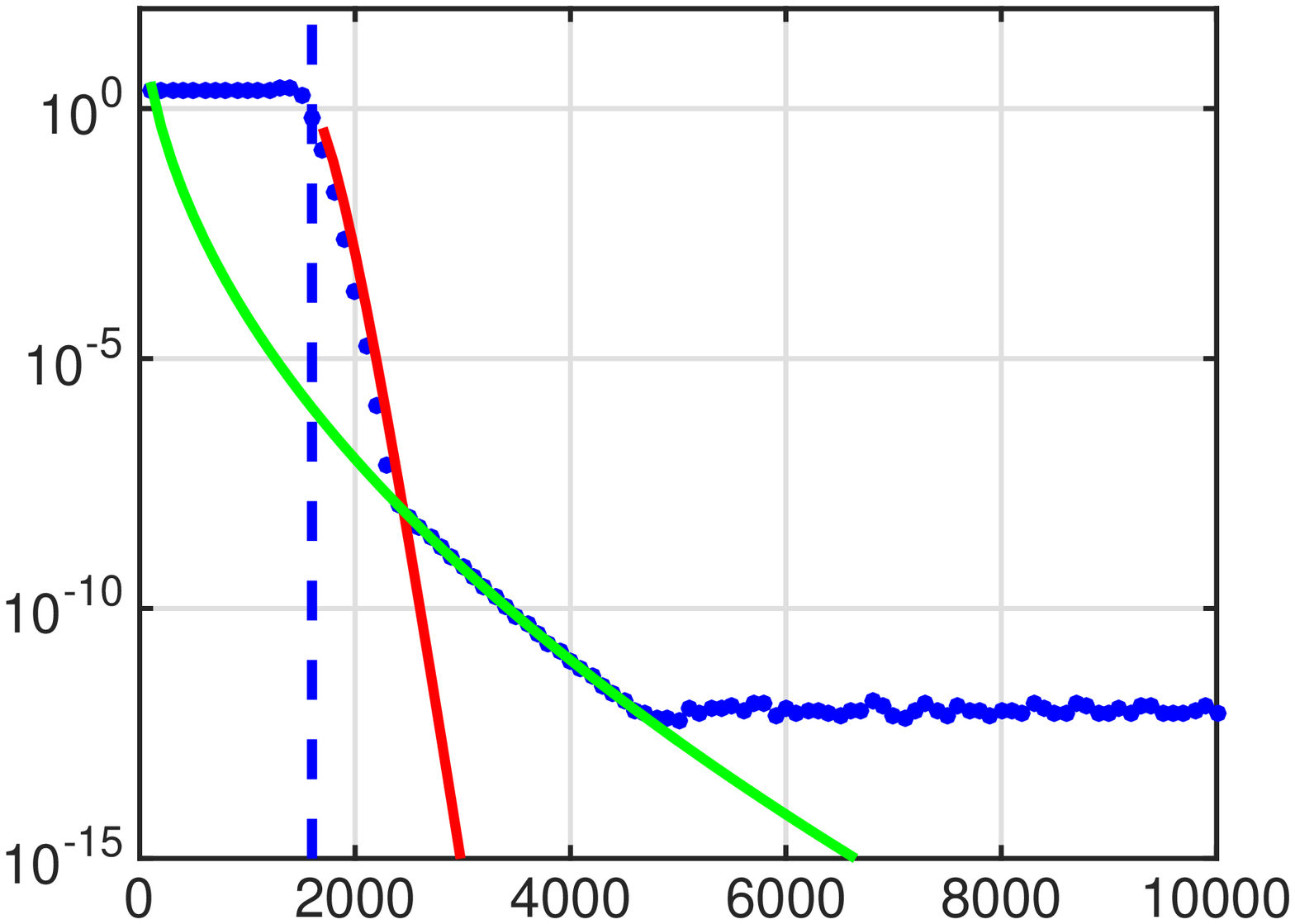}
&
\includegraphics[width=4.0cm]{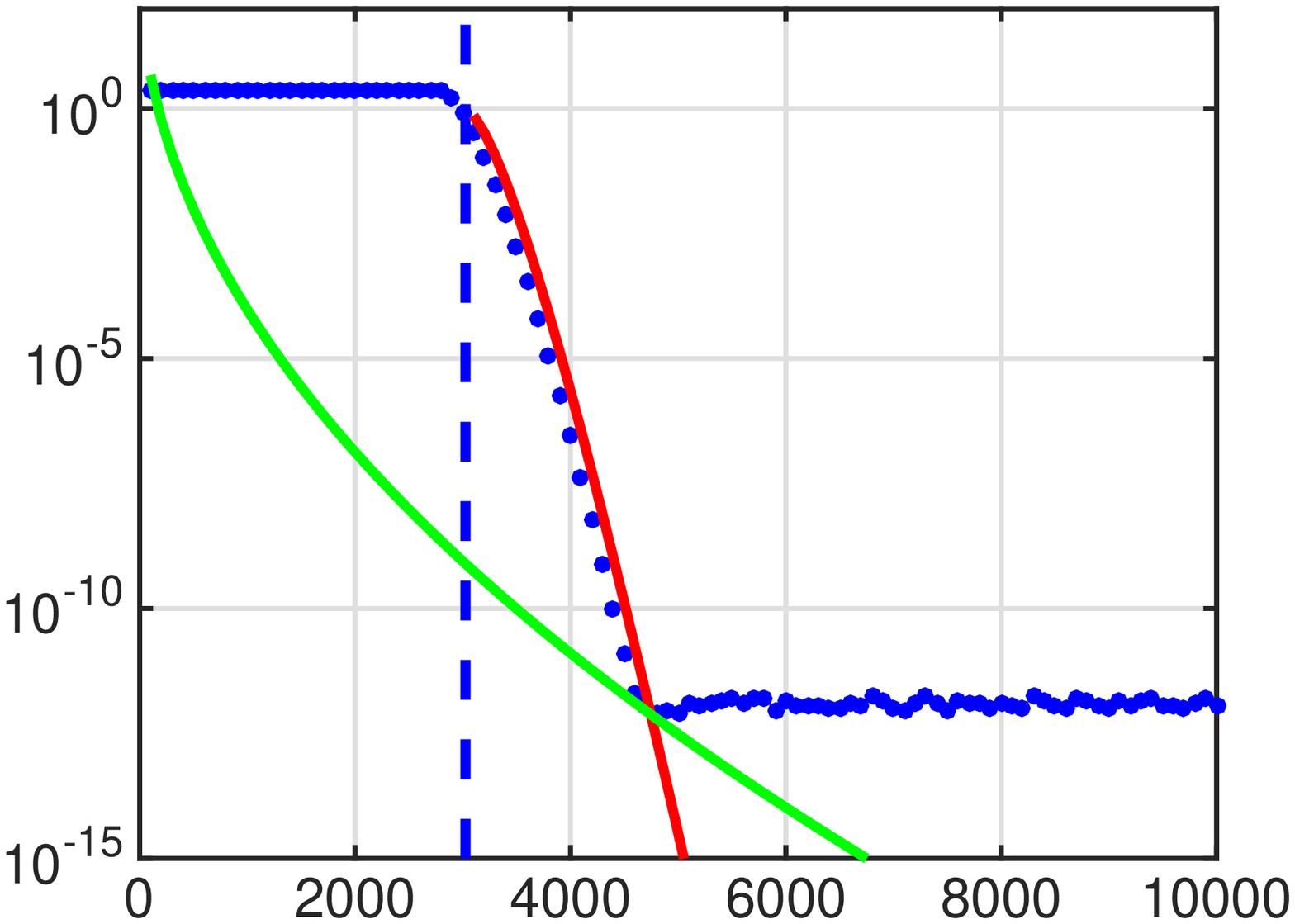}
&
\includegraphics[width=4.0cm]{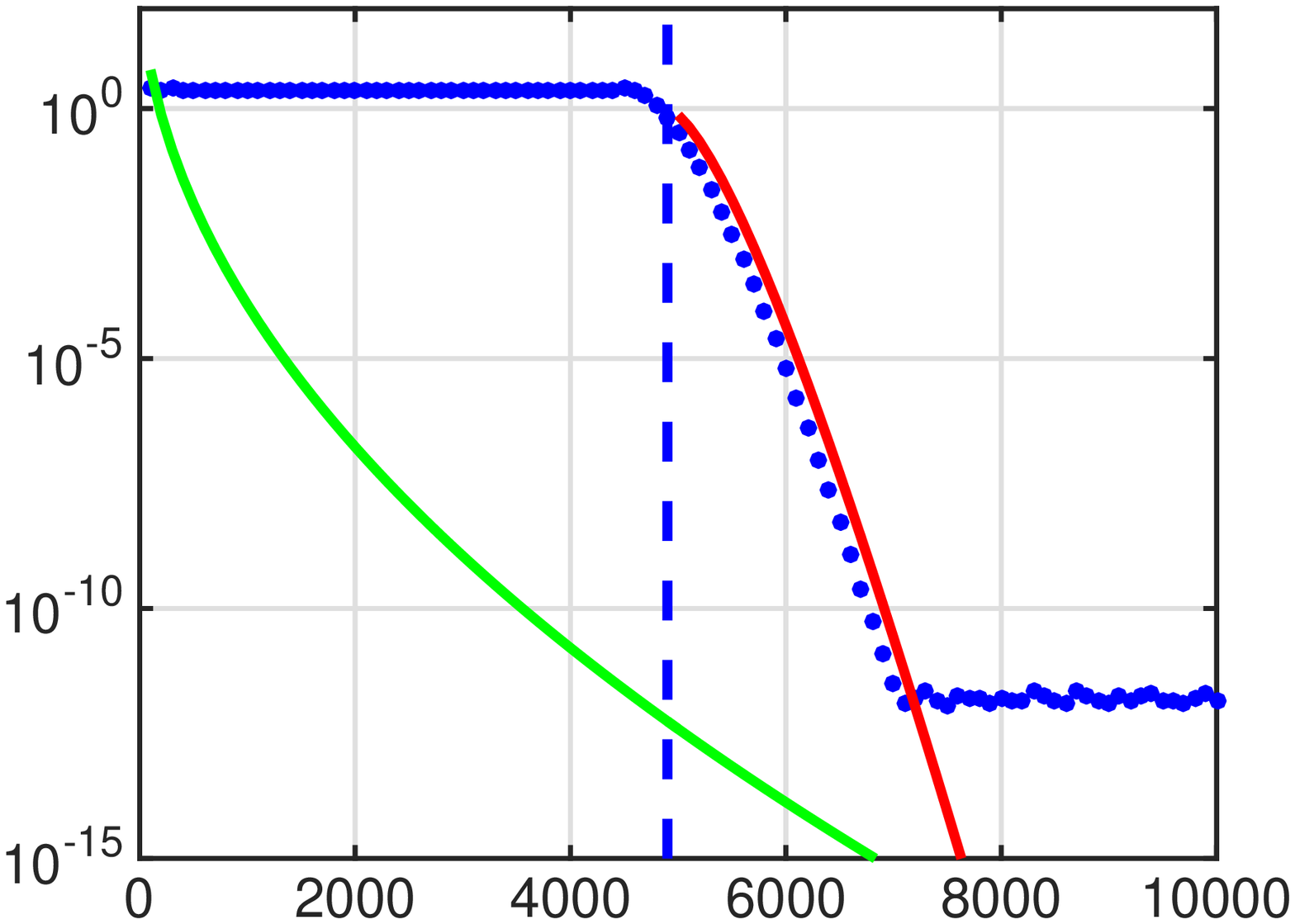}
\\
$\omega = 80$ & $\omega = 110$ & $\omega = 140$

\end{tabular}
\end{center}
\caption{The error against $n$ for $\psi_E$ with $L = c \sqrt{n}$, $c=0.5$, and the oscillatory function $f(x) = \exp(-2 \pi \I \omega x)$.  The dashed line is the theoretical resolution power $c^2$.  The red curve is the interior error bound $\exp(-\pi \omega H(\sqrt{n/n^*})$ and the green curve is the endpoint error bound $2 \pi \omega \exp(-c \sqrt{n})$.
}
\label{f:PsiEResConv}
\end{figure}

Theorem \ref{t:PsiE_respower} also allows us to examine much more precisely the behaviour of the error for oscillatory functions.  This is shown in Fig.\  \ref{f:PsiEResConv}.  Note the close agreement of the numerical error with the theoretical error bound \R{resn_err}.  One interesting facet of this diagram is a kink in the error that occurs for $\omega =80$.  This is due to the presence of the two exponentially-decaying terms $\exp(-\pi \omega H(\sqrt{n/n^*}))$ (which corresponds to the interior error) and $2 \pi \omega \exp(-c \sqrt{n})$ (corresponding to the endpoint error) in the error bound \R{resn_err}.  Since $H(t) \sim \pi t/2$ as $t \rightarrow \infty$,  the first term $\exp(-\pi \omega H(\sqrt{n/n^*})) \sim \exp(-\pi^2 \sqrt{n} / (2c))$ for large $n$.  Hence, for $c < \pi/\sqrt{2}$ this interior error term decays faster than the endpoint error term.  Thus, the interior error dominates for small $n$, but at a certain the endpoint error then begins to dominate, leading to the observed kink in the error graph.  Note that this kink may not be observed in finite-precision arithmetic, since the transition may occur when the error is already below machine epsilon.  This explains why it is not witnessed for larger values of $\omega$ in Fig.\ \ref{f:PsiEResConv}.

This kink phenomenon also explains why reducing $c$ in the exponential map so as to get the best resolution power -- see \R{E:PsiE_respower} -- is problematic.  Although the error does indeed begin to decay for a smaller value of $n$ when $c$ is small, if the kink phenomenon occurs then the error eventually decays at the slow rate of $\exp(-c \sqrt{n})$.

\section{The parametrized exponential map $\psi_{SE}$  }\label{s:psiSE}
We now turn our attention to the map \R{E:PsiSE}.

\subsection{Convergence and parameter choices}
For our main convergence result for this map, we first require the following lemma:

\lem{
\label{l:psiSEconv}
Let $\psi_{SE}$ be given by \R{E:PsiSE} for some $\alpha > 0$.  Suppose that $f$ is analytic and bounded in the region $\tilde{S}_{\alpha}$ given by \R{slitstrip}.  Let $M_{\psi}$ be as in \R{E:M1} with $\beta = \alpha$ and suppose that there exists $0 < \tau < 1$ such that \R{E:M2} holds with constant $N_{\psi}$.  Let $p_{n,L,\alpha}$ be as in \R{E:AproxPsi}.  If
\be{
\label{alpha_Lassumption}
0 < \alpha / (L-1/2) < 1,
}
then
\ea{
\| f - p_{n,L,\alpha} \|_{[0,1]} \leq    &  3 N_{\psi} \bar{\alpha}^{-1} \left ( \alpha \left | \log\left(1-\exp(\pi - \pi (L-1/2)/\alpha) \right) \right |  \right )^{\tau} \nn
\\
& + 114 M_{\psi} \bar{\alpha}^{-2} n \exp \left(-\alpha n \pi / (2L)  \right), \label{psiSEerr}
}
where $\bar{\alpha} = \min \{ \alpha , 1 \}$.
}
\prf{
We shall use Lemma \ref{lemConvPsiGeneral}.  Notice that
\bes{
C_{\psi} \leq \sup_{\substack{|x|,|y| \leq 1}} | \psi^{-1}_{SE}(-L+\alpha(x+\I y)) |,
}
and that
\bes{
\psi^{-1}_{SE}(-L+\alpha(x+\I y)) = \frac{\alpha}{\pi} \log \left ( \frac{1+\exp(\pi(1/2-L)/\alpha) \exp(\pi(x+\I y)) }{1+\exp(\pi(-1/2-L)/\alpha) \exp(\pi(x+\I y)) } \right ).
}
We now consider two cases:

\vspace{1pc}\noindent Case (i): \textit{$y = \pm 1$ fixed, $|x| \leq 1$ varying.}
If $y = \pm 1$ we have
\bes{
\psi^{-1}_{SE}(-L+\alpha(x\pm \I)) = \frac{\alpha}{\pi} \log \left ( \frac{1-\exp(\pi(1/2-L)/\alpha) \exp(\pi x) }{1-\exp(\pi(-1/2-L)/\alpha) \exp(\pi x) } \right ).
}
Note that
\bes{
0 < \frac{1-\exp(\pi(1/2-L)/\alpha) \exp(\pi x) }{1-\exp(\pi(-1/2-L)/\alpha) \exp(\pi x) }  < 1 ,
}
due to \R{alpha_Lassumption} and the fact that $|x| \leq 1$.  Also
\bes{
\frac{1-\exp(\pi(1/2-L)/\alpha) \exp(\pi x) }{1-\exp(\pi(-1/2-L)/\alpha) \exp(\pi x) } = \exp(\pi/\alpha) + \frac{1-\exp(\pi/\alpha)}{1-\exp(\pi(-1/2-L)/\alpha) \exp(\pi x)}.
}
Since $\alpha > 0$ this function is minimized at $x=1$ and therefore
\eas{
\sup_{|x| \leq 1} |\psi^{-1}_{SE}(-L+\alpha(x\pm \I)) | &= \frac{\alpha}{\pi} \left |\log \left (  \frac{1-\exp(\pi + \pi(1/2-L)/\alpha) }{1-\exp(\pi + \pi(-1/2-L)/\alpha)} \right ) \right |.
}
Note that the function $|\log(1-x)|$ is increasing on $0 < x < 1$.  Hence
\bes{
\sup_{|x| \leq 1} |\psi^{-1}_{SE}(-L+\alpha(x\pm \I)) | \leq \frac{2 \alpha}{\pi} \left | \log(1-\exp(\pi-\pi(L-1/2)/\alpha)) \right |.
}

\vspace{1pc}\noindent Case (ii): \textit{$x = \pm 1$ fixed, $|y| \leq 1$ varying.}
In this case, we have
\bes{
\psi^{-1}_{SE}(-L+\alpha(x\pm \I)) = \frac{\alpha}{\pi}  \log \left ( \frac{1+A \exp(\I \pi y) }{1+B \exp(\I \pi y) } \right ).
}
where $A = \exp(\pi(1/2-L)/\alpha \pm \pi)$ and $B = \exp(\pi(-1/2-L)/\alpha \pm \pi )$, and therefore
\bes{
\psi^{-1}_{SE}(-L+\alpha(\pm 1 + \I y)) \leq \frac{\alpha}{\pi} \left ( \left | \log ( 1+A \exp(\I \pi y)  ) \right | + \left | \log ( 1+B \exp(\I \pi y)  ) \right | \right ).
}
Consider $\log ( 1+A \exp(\I \pi y)  )$.
We have
\bes{
\log ( 1+A \exp(\I \pi y)  ) = \frac12 \log(1+2A \cos \pi y + A^2) + \I \tan^{-1} \left ( \frac{A \sin \pi y}{1+A \cos \pi y } \right ),
}
which gives
\bes{
\left | \log ( 1+A \exp(\I \pi y)  ) \right |^2 = \frac14 \left | \log(1+2A \cos \pi y + A^2) \right |^2 + \left ( \tan^{-1} \left ( \frac{A \sin \pi y}{1+A \cos \pi y } \right ) \right )^2.
}
For the first term, since $A < 1$ due to \R{alpha_Lassumption}, we deduce that
\bes{
\left | \log(1+2A \cos \pi y + A^2) \right | \leq 2 \left | \log(1-A) \right |
}
For the second term, we first note that $\tan^{-1}$ is an increasing function.  Hence we wish to maximize $\frac{A \sin \pi y}{1+A \cos \pi y }$.  After some simple calculus we deduce that this takes maximum value $\frac{A}{\sqrt{1-A^2}}$.  This gives
\bes{
\left | \log ( 1+A \exp(\I \pi y)  ) \right |^2 \leq | \log(1-A) |^2 + \left ( \tan^{-1} \left ( \frac{A}{\sqrt{1-A^2}} \right ) \right )^2.
}
To simplify matters, we now claim that
\bes{
| \log(1-x) | \geq \tan^{-1} \left ( x / \sqrt{1-x^2} \right ),\quad 0 \leq x < 1.
}
To see this, notice that both functions are zero when $x = 0$ and increasing in $x$.  Moreover, their derivatives are $\frac{1}{1-x}$ and $\frac{1}{\sqrt{1-x^2}}$ respectively.  Since $\frac{1}{1-x} \geq \frac{1}{\sqrt{1-x^2}}$, $0 \leq x < 1$, we have established the claim.  We now obtain
\bes{
\left | \log ( 1+A \exp(\I \pi y)  ) \right | \leq \sqrt{2} | \log(1-A) |.
}
For the term $\log(1+B\exp(\I \pi y))$ we follow an identical argument.  Noting that $| \log(1-x) |$ is an increasing function of $0 < x <1$, we now conclude that
\bes{
\sup_{|y| \leq 1 } \left | \psi^{-1}_{SE}(-L+\alpha(\pm1+ \I y)) \right | \leq \frac{2\sqrt{2}\alpha}{\pi}  | \log(1-\exp(\pi - \pi (L-1/2)/\alpha)) | .
}
Combining this with case (i), we deduce that
\bes{
(C_{\psi})^{\tau} \leq \left ( \alpha | \log ( 1 - \exp(\pi - \pi(L-1/2)/\alpha) ) | \right )^{\tau},
}
and the result now follows from Lemma \ref{lemConvPsiGeneral}.
%
}

Much like the standard exponential map, the bound \R{psiSEerr} allows us to determine choices for the parameters $L$ and $\alpha$.  If $L$ is uniformly bounded in $n$ and $\alpha / (L-1/2) \rightarrow 0$ as $n \rightarrow \infty$, then one readily deduces that the optimal parameter choices are
\bes{
\alpha = \alpha_0 / \sqrt{n},\qquad L = 1/2 + L_0,
}
for constants $\alpha_0,L_0 > 0$.  From this we obtain the following, which is our main result:

\thm{
\label{t:psiSEconv}
Let $f$ be analytic and bounded in $\tilde{S}_{\gamma}$ for some $\gamma > 0$.  Let $\alpha_0 , L_0 > 0$ be fixed and suppose that $\psi_{SE}$ is the mapping given by \R{E:PsiSE} with  $L = 1/2 + L_0$ and $\alpha = \alpha_0 / \sqrt{n}$.  Let $M_{\psi}$ be as in \R{E:M1} with $\beta = \alpha$ and suppose that there exists $0 < \tau < 1$ such that \R{E:M2} holds with constant $N_{\psi}$.  Then, for all $n \geq n_0 = \max \{ \alpha_0 , \alpha_0 / \gamma , 2 \alpha_0/L_0 \}^2$ we have
\be{
\label{psiSEerr2}
\| f - p_{n,L,\alpha} \|_{[0,1]} \leq A \left ( \frac{M_{\psi} n^2}{\alpha^2_0}  \exp(-\alpha_0 \pi \sqrt{n} / (2L_0+1) ) + \frac{N_{\psi} \sqrt{n}}{\alpha_0}  \exp ( -\tau \pi L_0 \sqrt{n} / \alpha_0 ) \right ),\qquad
}
where $M_{\psi}$ and $N_{\psi}$ are as in Lemma \ref{l:psiSEconv} and the constant $A \leq 114$.
}
\prf{
Observe that $L>1/2$ for all $n$, and $\alpha / (L-1/2) \leq 1/2$ and $\alpha \leq 1$ for $n \geq n_0$, and that $f$ is analytic in $\tilde{S}_{\alpha}$ for $n \geq n_0$.  Hence we may apply Lemma \ref{l:psiSEconv}, to get
\eas{
\| f - p_{n,L,\alpha} \|_{[0,1]} \leq & 3 N_{\psi} \alpha^{-1} \left | \log\left(1-\exp(\pi - \pi (L-1/2)/\alpha) \right) \right |^{\tau} 
\\
&+ 114 M_{\psi} \alpha^{-2} n \exp(-\alpha n \pi / (2 L)). 
}
Note that $| \log(1-x) | \leq \frac{| \log(1-a) |}{a} x$ for $0 \leq x \leq a <1$.  Hence, since $\alpha / (L-1/2) \leq 1/2$ we have
\ea{
| \log(1-\exp(\pi - \pi (L-1/2)/\alpha)) | & \leq \frac{| \log(1-\exp(-\pi))) |}{\exp(-\pi)} \exp(\pi - \pi (L-1/2)/\alpha)) \nn
\\
& \leq 24 \exp(-\pi (L-1/2)/\alpha). \label{logexpstep}
}
This gives
\bes{
\| f - p_{n,L,\alpha} \|_{[0,1]} \leq 72 N_{\psi} \alpha^{-1} \exp(-\tau \pi (L-1/2) / \alpha) + 114 M_{\psi} \alpha^{-2} n \exp(-\alpha n \pi / (2 L)).
}
We now substitute the values of $\alpha$ and $L$ to deduce the result.
}

Theorem \ref{t:psiSEconv} predicts root-exponential decay of the error with index
\be{
\label{rho_psiSE}
\rho = \min \left \{ \alpha_0 \pi / (2 L_0 + 1) , \tau \pi L_0 / \alpha_0 \right \}.
}
This in good agreement with the example shown in the right panel of Fig.\ \ref{f:PsiEConv}.  We note also that Remark \ref{r:pi_over_two} carries over to this setting as well (see \cite[Thm.\ 3.4]{AdcockRichardsonMappings} for the case of the Chebyshev-based method).

\rem{
If $\tau > 0$ is known, then one may optimize the parameters by equating the terms in \R{rho_psiSE} to give $\alpha_0 = \sqrt{\tau L_0 (2L_0+1)/\pi}$ and
\bes{
\| f - p_{n,L,\alpha} \| = \cO_{a} \left ( \exp\left(-\sqrt{\tau L_0 \pi / (2L_0+1)}\sqrt{n}\right) \right ),\quad n \rightarrow \infty.
}
At first sight, it therefore appears advisable to make $L_0$ large to obtain the fastest index of convergence.  However, as we show next, this will lead to worse resolution properties for oscillatory functions.  Moreover, in general taking $L_0$ large means that the asymptotic root-exponential decay of the error will take longer to onset, since the parameter $n_0$ in Theorem \ref{t:psiSEconv} scales quadratically with $L_0$.
}

\subsection{Resolution power}
We first require the following lemma:

\lem{
Let $\psi_{SE}$ be the mapping given by \R{E:PsiSE} and suppose that $f(x) = \exp(-2 \pi \I \omega x)$ for some $\omega \geq 1$.  If  $p_{n,L,\alpha}$ is the approximation defined by \R{E:AproxPsi}, then
\bes{
\| f - p_{n,L,\alpha} \|_{[0,1]} \leq 144 \pi \E^2 \omega  \exp(-\pi(L-1/2)/\alpha) + 114 \frac{n}{\alpha^2} \exp\left ( 2 \pi \omega \alpha -\alpha n \pi / (2L) \right )
}
provided $0 < \alpha < 1$ and $\alpha / (L-1/2) < \min \left \{ 1/2, 1/(1+\pi^{-1} \log(4 \alpha \omega)))\right \}$.
}
\prf{
As before, we commence by estimating the quantities $M_{\psi}$ and $N_{\psi}$ in Lemma \ref{l:psiSEconv} for the function $f(x) = \E^{-2 \pi \I \omega x}$ for large $ | \omega | $.  First, consider $M_{\psi}$.  By the maximum modulus principle, we have
\bes{
M_{\psi} = \sup_{z \in S_{\alpha}} | f(z)| = \sup_{\substack{z = x \pm \I \alpha \\ x \in \bbR }} \left| \exp \left( - 2\pi  \I \omega  (x+ \I \alpha) \right) \right|  = \exp \left( 2 \pi \omega \alpha \right)  }
Now let us study $ N_{\psi} $.  Since we may take $ \tau = 1$, we find that
\bes{
N_{\psi} \leq 2 \pi \omega \sup_{|z \pm L| \leq \alpha} | f(\psi^{-1}(z)) | = 2 \pi \omega \exp \left (2 \pi \omega \sup_{0 \leq \theta < 2 \pi} \Im \psi^{-1}\left (-L+\alpha \E^{\I \theta} \right ) \right ).
}
Note that
\bes{
\psi^{-1} ( - L + \alpha \E^{ \I \theta } ) = \frac{ \alpha } { \pi } \log \left( \frac{  1 +  \exp \left( \pi \left( 1/2 - L \right) / \alpha \right) \exp \left(   \pi \E^{\I \theta}   \right)  }  { 1 +   \exp \left( \pi \left( -1/2 - L \right) / \alpha \right) \exp \left(   \pi \E^{\I \theta }   \right)  } \right).
}
Consider the expression $\log(1+x+\I y)$, where $|x|,|y| \leq \delta$ and $0 < \delta < 1/2$.  We have
\bes{
\left | \Im \log(1+x+\I y) \right | =\left | \tan^{-1} \left ( \frac{y}{1+x} \right ) \right | \leq \tan^{-1} \left ( \frac{\delta}{1-\delta} \right ) \leq 2 \delta.
}
Now set $x+\I y = \exp \left( \pi \left( \pm 1/2 - L \right) / \alpha \right) \exp \left(   \pi \E^{\I \theta}  \right )$.  Then we may take
\bes{
\delta = \exp(\pi + \pi(1/2-L)/\alpha) < 1/2,
}
and it follows that
\bes{
\Im \psi^{-1}(-L+\alpha \E^{\I \theta}) \leq 4 \alpha / \pi \exp(\pi -\pi(L-1/2)/\alpha).
}
This now gives
\bes{
N_{\psi} \leq 2 \pi \omega \exp \left ( 8 \alpha \omega \exp(\pi -\pi(L-1/2)/\alpha) \right ) \leq 2 \pi \E^{2} \omega,
}
due to the assumptions on $\alpha / (L-1/2)$.  We now apply Lemma \ref{l:psiSEconv} and \R{logexpstep}.
}

From this we deduce our main result:

\thm{ \label{t:PsiSE_respower}
Let $\psi_{SE}$ be the mapping given by \R{E:PsiSE} with $\alpha = \alpha_0/\sqrt{n}$ and $L = L_0+1/2$ for $\alpha_0,L_0 > 0$ and suppose that $f(x) = \exp(-2 \pi \I \omega x)$ for some $\omega \geq 1$.  If  $p_{n,L,\alpha}$ is  the approximation defined by \R{E:AproxPsi}, then
\be{
\label{resn_err_PsiSE}
\| f - p_{n,L,\alpha} \|_{[0,1]} = \cO_a \left ( \exp(\alpha_0 \pi(2 \omega - n / (2L_0+1))/\sqrt{n})+ \omega \exp\left (-\pi L_0 \sqrt{n} / \alpha_0 \right ) \right ),
}
as $n \rightarrow \infty$, uniformly in $\omega$.  In particular, the resolution power satisfies
\be{
\label{PsiSE_respower}
\limsup_{\omega \rightarrow \infty} \cR(\omega ; \delta) / \omega \leq 4(L_0+1/2),\qquad 0 < \delta < 1.
}
}

\begin{figure}
\begin{center}
\begin{tabular}{ccc}
\includegraphics[width=4.0cm]{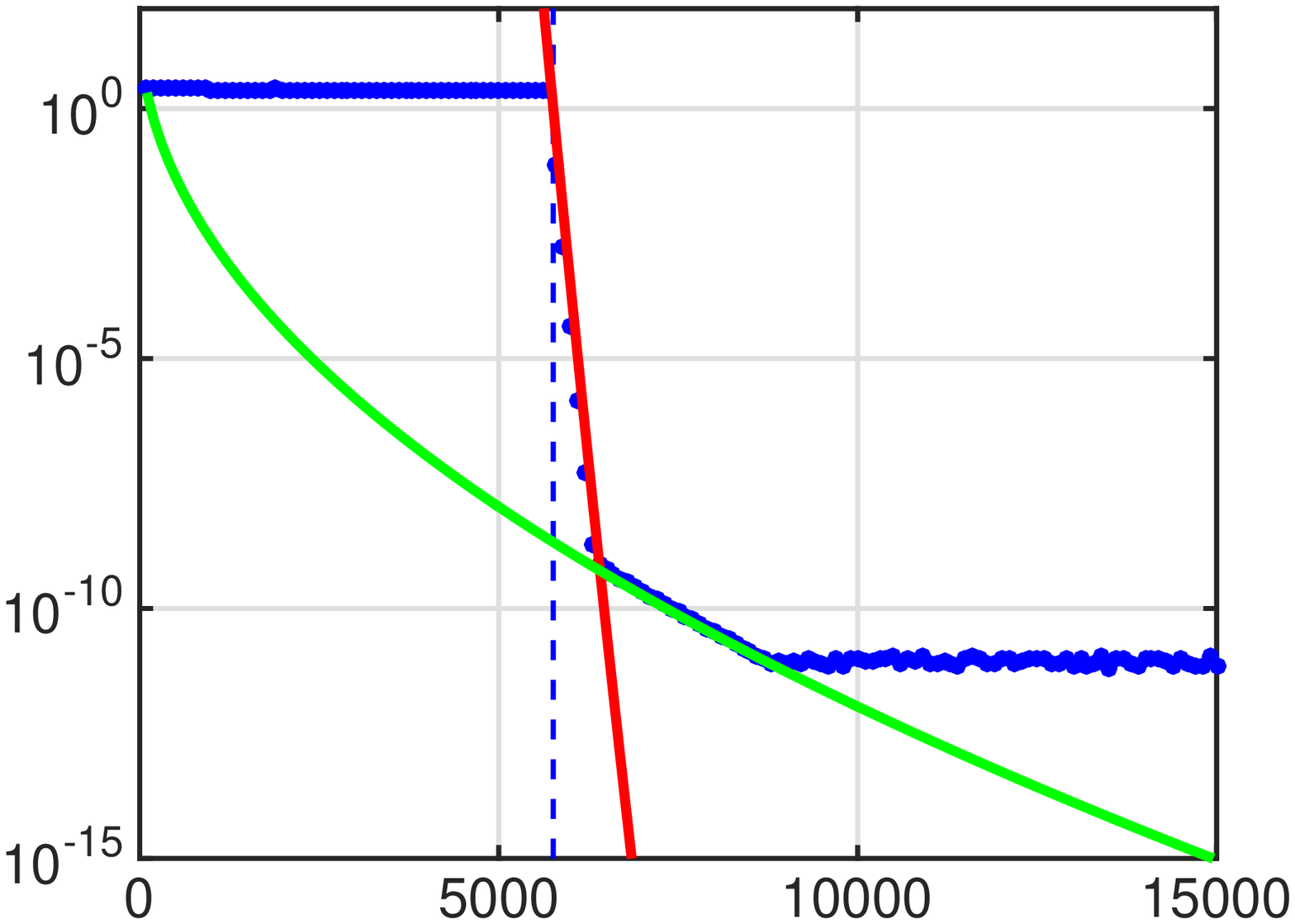}
&
\includegraphics[width=4.0cm]{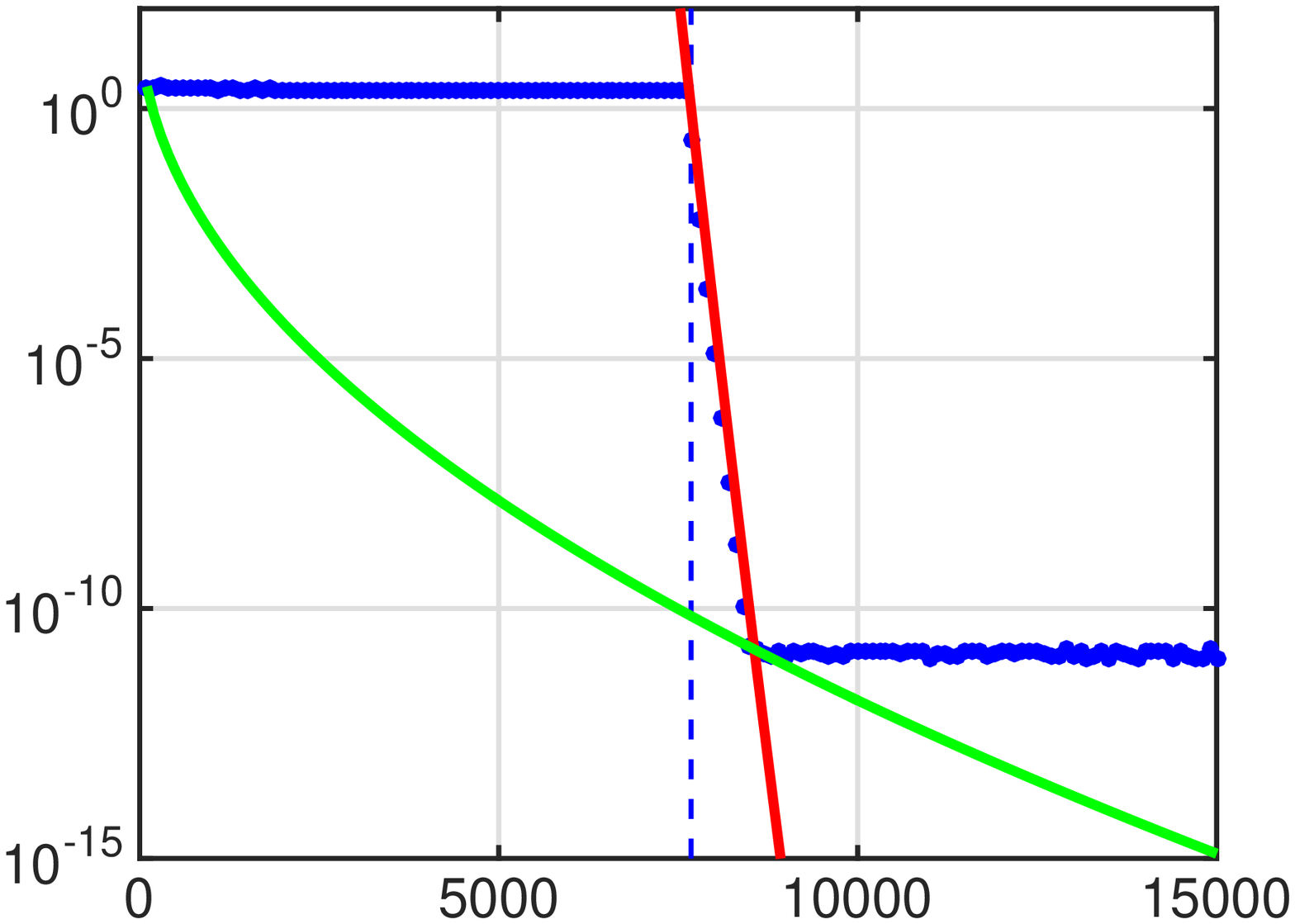}
&
\includegraphics[width=4.0cm]{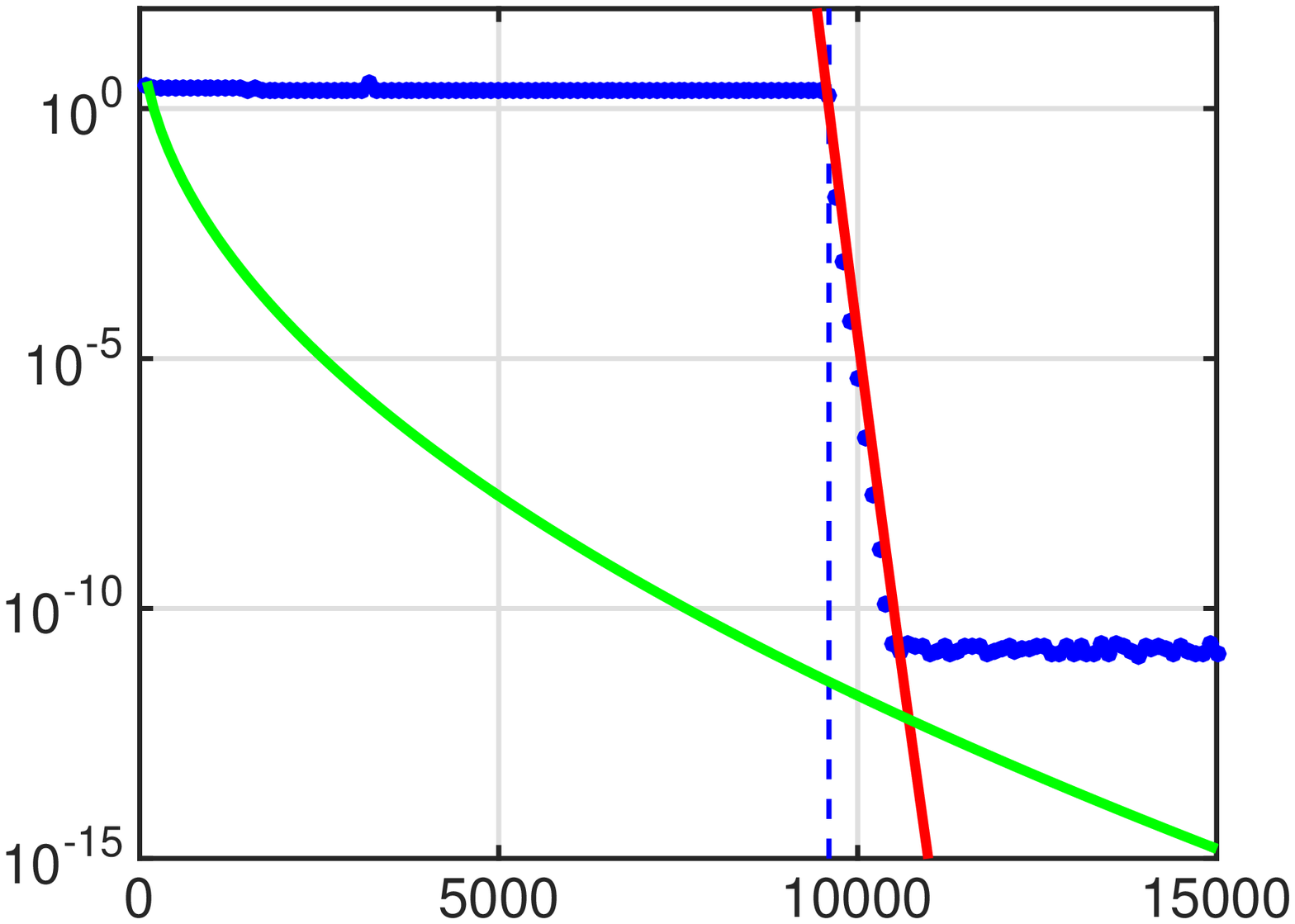}
\\
$\omega = 2400$ & $\omega = 3200$ & $\omega = 4000$

\end{tabular}
\end{center}
\caption{The error against $n$ for $\psi_{SE}$ with $L_0 = 0.1$, $\alpha_0 = 1$ and the oscillatory function $f(x) = \exp(-2 \pi \I \omega x)$.  The dashed line is the theoretical resolution power $4(L_0+1/2)$.  The red and green curves are the interior and endpoint error bounds respectively, given by \R{resn_err_PsiSE}.
}
\label{f:PsiSEResConv}
\end{figure}

In the right panel of Fig.\  \ref{f:PsiEResolution} we show the numerical verification of the resolution power \R{PsiSE_respower}.  Theorem \ref{t:PsiSE_respower} also allows us to understand the behaviour of the error for oscillatory functions in more detail.  This is shown in Fig.\  \ref{f:PsiSEResConv}.  The numerical error is in good agreement with the theoretical error bound \R{resn_err}.
Again, we witness a kink phenomenon for a range of $ \omega $ values.  We note also the clear superiority of this method over $\psi_E$ for oscillatory functions (compare with Fig.\ \ref{f:PsiEResConv}).

\section{The double-exponential map $\psi_{DE}$}\label{s:DE}

In this and \S \ref{s:psiSDE} we consider double-exponential maps.  The convergence and resolution analysis is more difficult for these maps, making precise statements (valid for any $n \in \bbN$) harder to obtain.  Hence, we now opt for a less formal approach in which we derive asymptotic error bounds as $n \rightarrow \infty$.  Note, however, that the corresponding results on resolution power remain sharp even with this approach.

We first consider convergence of the double-exponential map $\psi_{DE}$.  For this, we make use of the Lambert-W function \cite{CorlessEtAlLambertW}.  Recall that $W(x)$ is defined implicitly by the relation $x = W(x) \exp(W(x))$ and on its principal branch satisfies $W(x) \sim \log x - \log \log x$ as $x \rightarrow \infty$.

\thm{
\label{t:PsiDEConv}
Let $\psi_{DE}$ be the mapping given by \R{E:PsiDE}. Suppose that $f(x)$ is analytic and bounded in the domain $\psi^{-1}(S_{\beta})$ for some $0 < \beta < 1$ and let $M_{\psi}$ and $N_{\psi}$ be as in \R{E:M1} and \R{E:M2} respectively.  Let $p_{n,L}$ be the approximation defined in \R{E:AproxPsi}, where $L = 1+W(cn)$ for some $c>0$.  Then
\be{
 \label{E:PsiDEConv}
\| f - p_{n,L}\|_{[0,1]} = \cO_a \left ( M_{\psi} \rho^{-n/\log(cn)}_1+N_{\psi} \rho^{-n/\log(cn)}_2 \right ),\qquad n \rightarrow \infty,
}
where $\rho_1 = \exp(\beta \pi/2)$ and $\rho_2 = \exp(\tau \pi c /2)$.
}

\prf{
Notice that $\psi^{-1}$ is analytic in $S_{\beta}$ for any $\beta < \pi/2$.  We now apply Lemma \ref{lemConvPsiGeneral}.  The first term in \R{E:PsiDEConv} follows immediately from the corresponding term in \R{E:lemma}.  For the second, we need to estimate $C_{\psi}$.  For this, observe that
\eas{
\psi^{-1}(-L+\beta \E^{\I \theta} ) & \sim \left ( 1 + \exp \left ( \pi \exp(L-\beta \E^{\I \theta}) / 2 \right ) \right )^{-1}
\sim  \exp \left ( -\pi \exp(L-\beta \E^{\I \theta}) / 2 \right ),
}
as $n \rightarrow \infty$, where in the second step we use the fact that $\beta < 1$, so that $\Re (\exp(L-\beta \E^{\I \theta})) = \exp(L-\beta \cos \theta) \cos (\beta \sin(\theta)) \geq \exp(L-\beta) \cos(\beta) > 0$ is strictly positive for all $0 \leq \theta < 2 \pi$.  Hence, as $n \rightarrow \infty$,
\eas{
C_{\psi} & \sim \max_{0 \leq \theta < 2\pi} \left | \exp \left ( -\pi \exp(L-\beta \E^{\I \theta}) / 2 \right ) \right |
\\
&= \exp \left ( -\pi \E^{L}/2 \min_{0 \leq \theta < 2 \pi} \exp(-\beta \cos(\theta)) \cos(\beta \sin(\theta)) \right )
}
It is readily checked that the function $g(\theta) = \cos(\beta \sin(\theta)) \E^{-\beta \cos(\theta)}$ satisfies $g'(\theta) = 0$ if and only if $\sin(\theta-\beta \sin(\theta)) = 0$.  If $0 < \beta < 1$ then $\min g(\theta) = g(0) = \exp(-\beta) \geq \exp(-1)$.  Hence we now get
\bes{
C_{\psi} \sim \exp(-\pi \E^{L-1}/2 ) \sim \exp(-c n \pi  / (2 \log(cn))),\qquad n \rightarrow \infty,
}
as required.
}

Next we address resolution power:

\thm{
  \label{t:PsiDE_respower}
Let $\psi_{DE}$ be the mapping given by \R{E:PsiDE} and let $p_{n,L}$ the approximation defined by \R{E:AproxPsi}.  If $L = 1+W(cn)$ then the resolution power satisfies
\be{
\label{PsiDE_respower}
\limsup_{\omega \rightarrow \infty} \cR(\omega) / \left ( \omega \log(c \omega)\right ) \leq \pi .
}
}
\prf{
We use the error estimate \R{E:PsiDEConv}.  Since this is valid for any $0 < \beta < 1$ we shall consider the asymptotic regime $\beta \rightarrow 0$.  In the usual manner,
\bes{
M_{\psi} = \exp \left ( 2 \pi \omega \sup_{\substack{z = x \pm \I y \\ x \in \bbR, |y| \leq \beta}} \Im \left ( \frac{1}{1 +\exp(-\pi \sinh(z))} \right ) \right ).
}
Let $z = x\pm\I y$ and write $ \Im \left ( \frac{1}{1 +\exp(-\pi \sinh(z))} \right ) = \pm g(x,y)$, where
\bes{
g(x,y) =  \frac{ \sin( \pi \cosh(x)  \sin( y) ) }  { 2 ( \cos( \pi \cosh(x)  \sin( y) ) + \cosh( \pi \sinh(x)  \cos( y) )   ) } .
}
Since $g(x,y)$ is an even function of $x$ we need only consider $x \geq 0$.  Also, note that
\bes{
|g(x,y)| \leq \frac{1}{2(\cosh(\pi \sinh(x) \cos(1))-1)},
}
Hence $|g(x,y) | \leq \beta \pi / 4$ for all $|y| \leq \beta$ provided $x \leq x_\beta$, where
\be{
\label{x_beta}
x_{\beta} =\sinh^{-1} \left ( \cosh^{-1} \left ( 1+ 2/(\beta \pi) \right ) / \left ( \pi \cos(1) \right ) \right ) \sim \log(\log(\beta)),\qquad \beta \rightarrow 0.
}
Now consider the behaviour of $g(x,y)$ for $0 \leq x \leq x_{\beta}$.  As $\beta \rightarrow 0$, we have
\bes{
g(x,y) \sim \frac{\pi \cosh(x) y}{2 \left ( 1 + \cosh(\pi \sinh(x)) \right ) } + \ord{(\beta \log \beta)^2},
}
and this holds uniformly in $|y| \leq \beta$ and $0 \leq x \leq x_{\beta}$ due to \R{x_beta}.  The function $h(x) = \frac{\cosh(x) }{1 + \cosh(\pi \sinh(x))}$ satisfies $h(0) = 1/2$ and $h(x) \rightarrow 0$ as $x \rightarrow \infty$.  Also, for $x > 0$,
\eas{
h'(x) &= \frac{\cosh(x) \left ( \tanh(x) - \pi \cosh(x) \tanh (\pi \sinh(x) / 2) \right ) }{1+\cosh(\pi \sinh(x))} < 0
}
since $\tanh(t)$ is an increasing function.  Hence $h$ attains its maximum value at $x=0$.  Combining this with the previous estimates, we deduce that
\bes{
\sup_{x \geq 0, | y | \leq \beta} g(x,y) \sim \beta \pi / 4,\qquad \beta \rightarrow 0,
}
and this gives $M_{\psi} \sim \exp ( \pi^2 \omega \beta / 2)$ as $\beta \rightarrow 0$.

We now consider $N_{\psi}$.  We have
\bes{
N_{\psi} \leq 2 \pi \omega \exp \left ( 2 \pi \omega \sup_{0 \leq \theta < 2 \pi} \Im \left ( \frac{1}{1+\exp(-\pi\sinh(-L+\beta \E^{\I \theta}))} \right ) \right ).
}
Now
\eas{
 \Im &\left ( \frac{1}{1+\exp(-\pi\sinh(-L+\beta \E^{\I \theta}))} \right )
 \\
 & = \frac{  \sin(\pi \cosh(  L - \beta \cos \theta )  )  \sin ( \beta \cos \theta )  }   {  2 (  \cos( \pi \cosh ( L - \beta \cos ( \theta )  )  \sin (  \beta \sin( \theta )  )   )      + \cosh( \pi \cos (  \beta \sin  \theta   )    )  \sinh (  L - \beta \cos  \theta )    ) }
 \\
 & \leq \frac{1}{2(\cosh(\pi \cos(1) \sinh(L-1)) - 1 ) }  \sim \frac{1}{2} \exp \left ( - c \pi \cos(1) n / 2 \right ),
}
as $n \rightarrow \infty$.
Combining this with  Theorem \ref{t:PsiDEConv} and the estimate for $M_{\psi}$ we get
\eas{
\| f - p_{n,L} \|_{[0,1]} = \cO_a \Big ( & \exp \left ( \beta \pi / 2 ( \pi \omega - n / \log(cn) ) \right )
\\
& +\omega \exp \left ( 2 \pi \omega \exp(-c \pi \cos(1) n / 2) -  \pi c n /(2 \log(cn)) \right ) \Big ),
}
as $n \rightarrow \infty$ and $\beta \rightarrow 0$.  Observe that
\bes{
2 \pi \omega \exp \left ( 2 \pi \omega \exp(-c \pi \cos(1) n / 2) \right ) \lesssim 1,\qquad n \gtrsim \frac{2}{c \pi \cos(1)} \log \left (2 \pi \omega / \log(2 \pi \omega) \right ).
}
Hence the second term in the error estimate is exponentially small once $n$ is logarithmically large in $\omega$.  Conversely, the first term $\exp \left ( \beta \pi / 2 ( \pi \omega - n / \log(cn) ) \right ) $ only begins to decay exponentially once $n/ \log(cn) \geq \pi \omega$, which is equivalent to $n \geq \pi \omega \log(c \omega)$.  This now gives the result.
}

\begin{figure}
\begin{center}
\begin{tabular}{ccc}
\includegraphics[width=5.5cm]{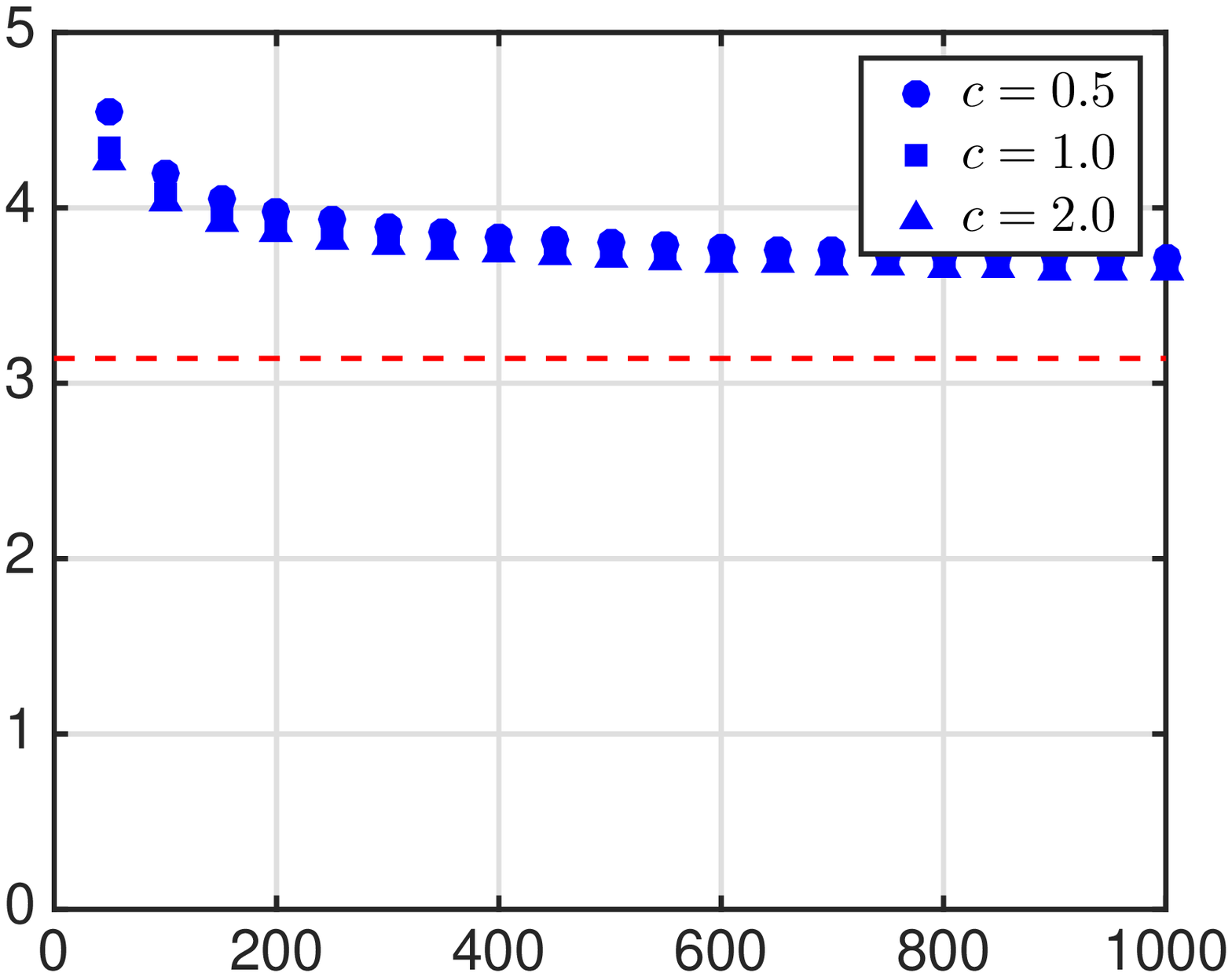} && \includegraphics[width=5.5cm]{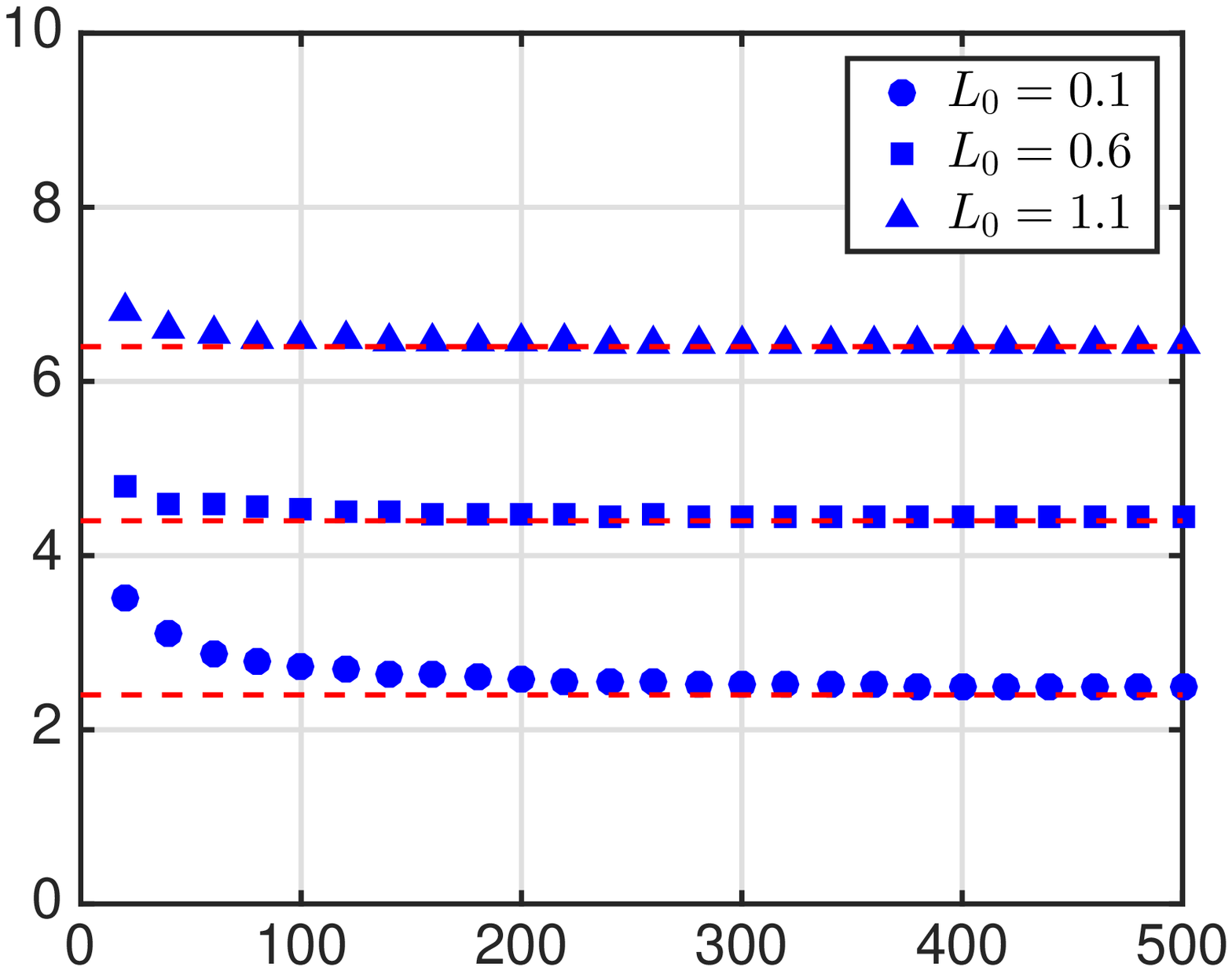}
\\
{\small $\psi_{DE}$, $L = W(cn)$} && {\small $\psi_{SDE}$, $L = L_0+1/2$, $\alpha = L_0 \pi / W(c n)$, $c=1$}
\end{tabular}
\end{center}
\caption{Left: the quantity $\cR(\omega ; \delta) / (\omega \log(c \omega))$ against $\omega$ with $\delta = 10^{-2}$ for $\psi_{DE}$.  The dashed red line shows the theoretical resolution constant $\pi$. Right: the quantity $\cR(\omega ; \delta) / \omega$ with $\delta = 10^{-2}$ for $\psi_{SDE}$.  The dashed red line shows the theoretical resolution constant $4(L_0+1/2)$.
}
\label{f:PsiDEResolution}
\end{figure}

Numerical verification of this theorem is given in Fig.\ \ref{f:PsiDEResolution}.

\section{The parametrized double exponential map $\psi_{SDE}$  }\label{s:psiSDE}

We commence with the following:

\lem{
\label{l:psiSDEconv}
For $\alpha > 0$ let $\psi_{SDE}$ be the mapping given by \R{E:PsiSDE}.  Suppose that $\alpha \rightarrow 0$ and that $L$ is uniformly bounded above and $L-1/2 >0$ is uniformly bounded away from zero.  Then the constant $C_{\psi}$ satisfies
\bes{
C_{\psi}(\sigma \alpha,L) \sim \frac{\alpha}{\pi}  \exp\left ( \pi/(2\alpha)-\exp(\pi(L-1/2)/\alpha-\sigma \pi) \right ),
}
uniformly in $0 < \sigma \leq \sigma_0 < 1/2$.
}
\prf{
The map $\psi^{-1}_{SDE}$ is analytic on the domain $S_{\beta}$ for any $\beta = \sigma \alpha$ with $\sigma < 1/2$ (see Lemma \ref{l:SDEanalytic}).  Now let $z = -L + \beta \E^{\I \theta}$.  Then, as $\alpha\rightarrow 0$, we have
\bes{
g^{-1}(z ; \alpha) \sim - \frac{\alpha}{\pi} \exp(\pi(L-1/2) / \alpha ) \exp(-\sigma \pi  \exp(\I \theta) ) = x+\I y,
}
where
\eas{
x &= - \frac{\alpha}{\pi} \exp(\pi (L-1/2) / \alpha) \exp(-\sigma \pi \cos(\theta) ) \cos(\sigma \pi  \sin(\theta) )
\\
y &= \frac{\alpha}{\pi} \exp(\pi (L-1/2) / \alpha ) \exp(-\sigma \pi  \cos(\theta)  ) \sin(\sigma \pi  \sin(\theta) ).
}
Note that $x$ is exponentially large as $\alpha \rightarrow 0$ and negative, since $\cos(\sigma \pi  \sin(\theta)) \geq \cos(\sigma_0 \pi) > 0$.  Hence,
\eas{
\psi^{-1}_{SDE}(z ; \alpha ) & = \frac{\alpha}{\pi} \log \left ( \frac{1+\exp(\pi(x+\I y + 1/2) / \alpha)}{1+\exp(\pi(x+\I y - 1/2) / \alpha)} \right ) \sim \frac{2 \alpha}{\pi} \exp(\pi(x+\I y ) / \alpha) \sinh(\pi / (2\alpha)),
}
from which it follows that
\bes{
C_{\psi} \sim \frac{\alpha}{\pi} \exp(\pi/(2\alpha)) \exp \left ( - \exp(\pi(L-1/2)/\alpha) \min_{0 \leq \theta < 2 \pi} \exp(-\sigma \pi  \cos(\theta) ) \cos(\sigma \pi  \sin(\theta)) \right ).
}
Recall from the proof of Theorem \ref{t:PsiDEConv} that $g(\theta) = \exp(-\sigma \cos(\theta)) \cos( \sigma \sin(\theta))$ is minimized at $\theta = 0$ since $0 < \sigma < 1/2$.  Hence we now get the result.
}

This leads to our main result on convergence:

\thm{
Let $\psi_{SDE}$ be the mapping defined by \R{E:PsiSDE} with $L = L_0 + 1/2$ and $\alpha = L_0 \pi / (\pi/2 + W(c n))$, where $L_0, c>0$ are fixed.  If $p_{n,L,\alpha}$ is the approximation defined by \R{E:AproxPsi} then
\bes{
\| f - p_{n,L,\alpha} \|_{[0,1]} = \cO_a \left ( (N_\psi + M_{\psi}) \rho^{-n / \log(cn)}  \right ),\qquad n \rightarrow \infty,
}
where $\rho = \max \left \{ \exp \left ( \pi^2 L_0 / (4 L_0 + 2) \right ) , \exp (c\tau) \right \}$.
}
\prf{
Substituting the parameter choices into Lemma \ref{l:psiSDEconv} we find that
\bes{
C_{\psi} = \cO \left (  \exp ( - c n /W(cn) ) \right ),\qquad n \rightarrow \infty,
}
for any $0 < \sigma < 1/2$.  We now apply Lemma \ref{lemConvPsiGeneral} to get that
\bes{
\| f - p_{n,L,\alpha} \|_{[0,1]} = \cO_a \left ( N_{\psi} \rho^{-n / \log(cn)}_1 + M_{\psi} \rho^{-n / \log(cn)}_2\right ),
}
where $\rho_1 = \exp(c \tau)$ and $\rho_2 = \exp(\sigma \pi^2 L_0 / (2L_0+1) )$.
}

We now turn our attention towards resolution power:

\thm{
Let $\psi_{SDE}$ be the mapping defined by \R{E:PsiSDE} with $L=L_0+1/2$ and $\alpha = L_0 \pi / (\pi/2 + W(cn))$ where $L_0,c>0$ are fixed.  Then the resolution power satisfies
\bes{
\limsup_{\omega \rightarrow \infty} \cR(\omega) / \omega \leq 4 L_0 + 2.
}
}
\prf{
By Lemmas \ref{l:psiSDEconv} and \ref{PsiSDEres_step} we have
\eas{
 M_{\psi}(\sigma \alpha) \exp(-\sigma \alpha n \pi / (2L) ) &\leq \exp(2 \pi \omega \sigma (  \alpha + o(\alpha)) - \sigma \alpha n \pi / (2L))
 \\
 & = \exp ( 2 \pi \sigma \alpha (\omega  - n / (4 L)) + o(\omega \sigma \alpha) ),
 }
 as $\alpha \rightarrow 0$ (i.e.\ $n \rightarrow \infty$) uniformly in $0 < \sigma \leq \sigma^*$.  Set $\sigma = 1/\omega$.  Then
 \bes{
 M_{\psi}(\sigma \alpha) \exp(-\sigma \alpha n \pi / (2L) ) \leq \exp(2 \pi \alpha(1 - n/(4 L \omega) ) + o(\alpha) ),\qquad \alpha \rightarrow 0.
 }
 Hence, the right-hand side begins to decay once $n \geq 4 L \omega$.  Also, we have
 \bes{
 N_{\psi} (\sigma \alpha,L) C_{\psi}(\sigma \alpha , L) \leq 2 \alpha \omega (1+o(1)) \exp(\pi/(2 \alpha) - \exp(\pi(L-1/2)/\alpha - \sigma \pi))
 }
Using the values for $L$ and $\alpha$, we see that this term is negligible as soon as $n$ is on the order of $\log \omega$.  Hence the result follows.
}

Numerical verification of this theorem is given in Fig.\ \ref{f:PsiDEResolution}.

\section{Numerical comparisons}\label{s:numerical}
\begin{figure}
\begin{center}
\begin{tabular}{ccc}
\includegraphics[width=4cm]{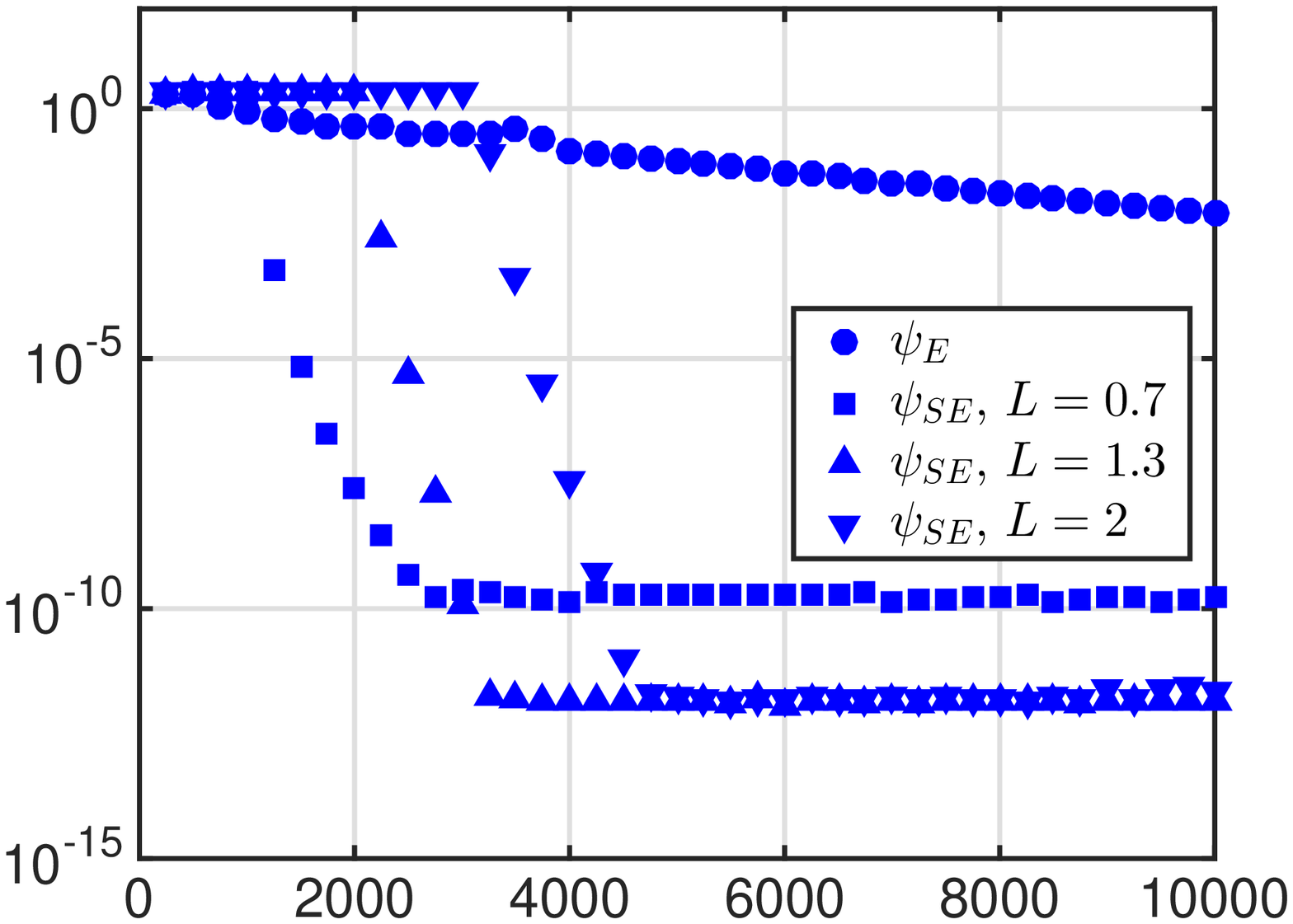}
&
\includegraphics[width=4cm]{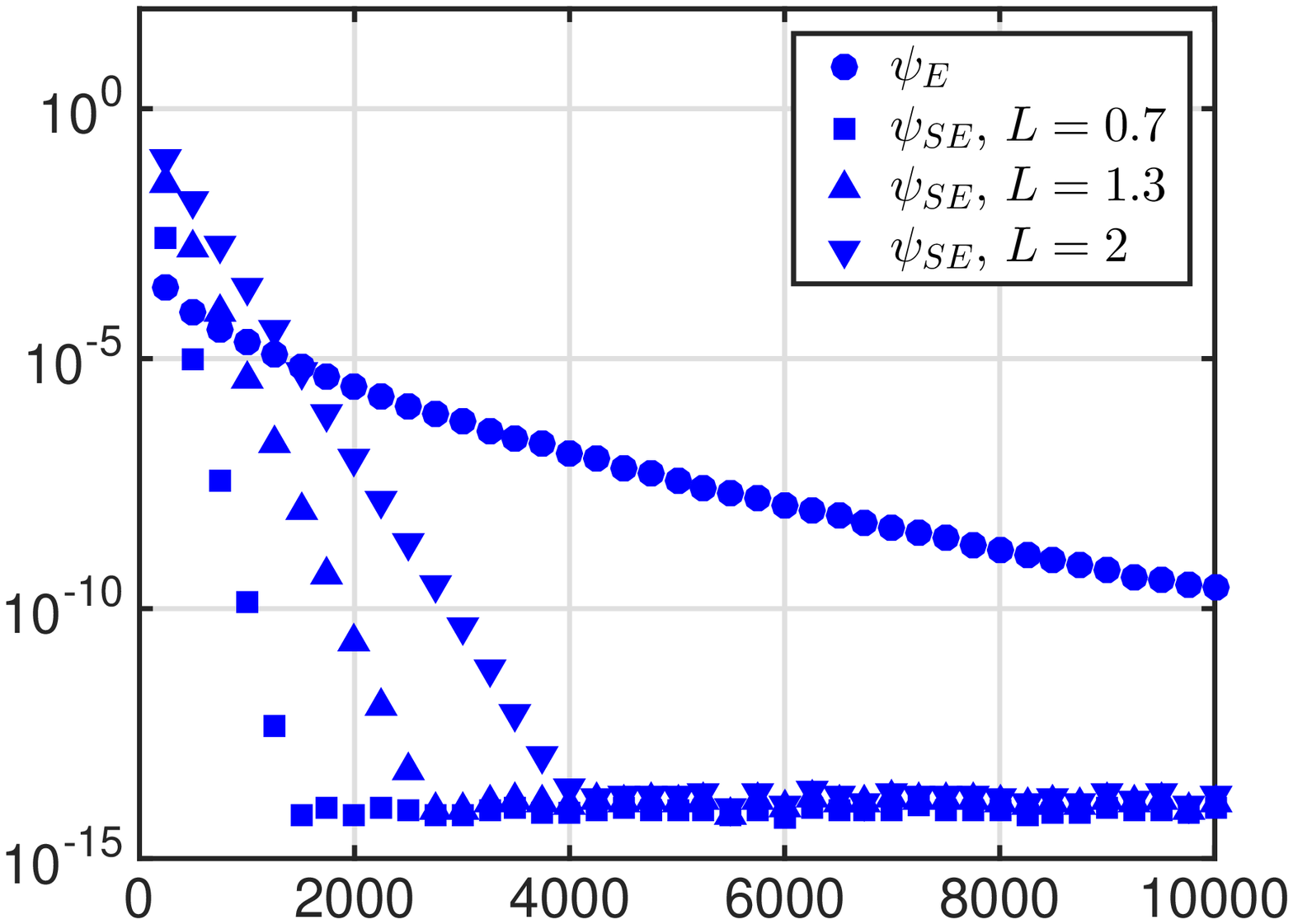}
&
\includegraphics[width=4cm]{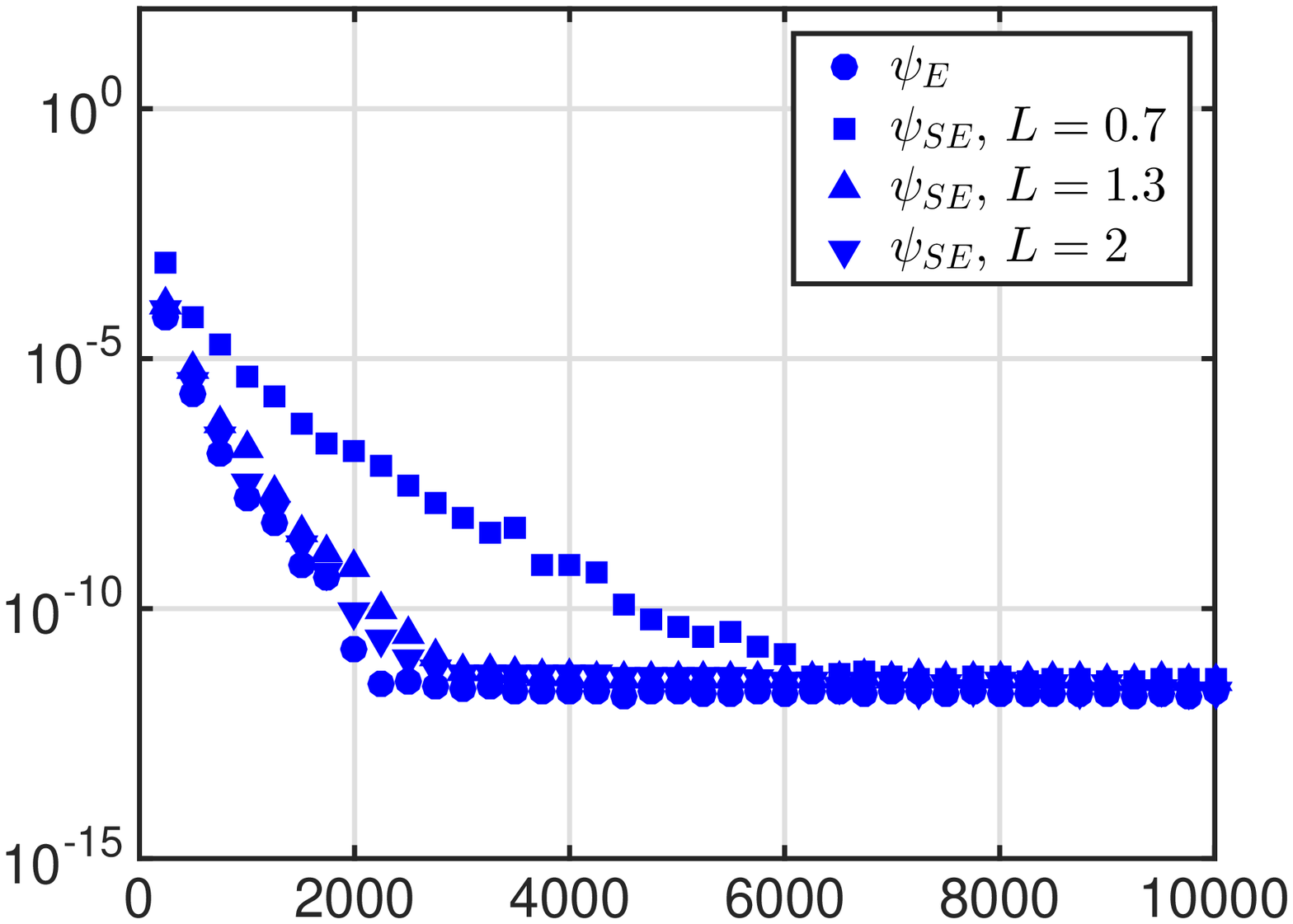}
\\
\includegraphics[width=4cm]{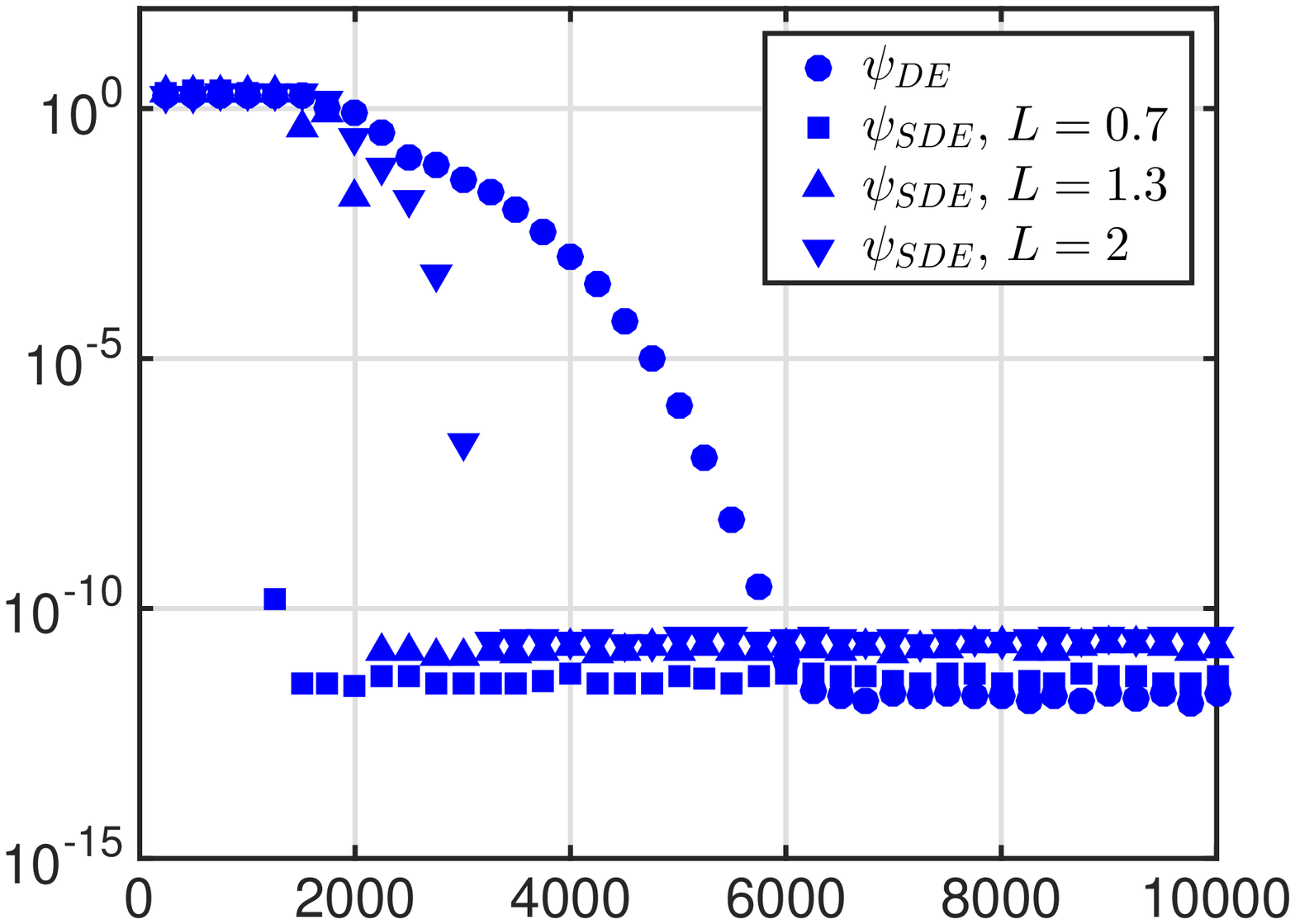}
&
\includegraphics[width=4cm]{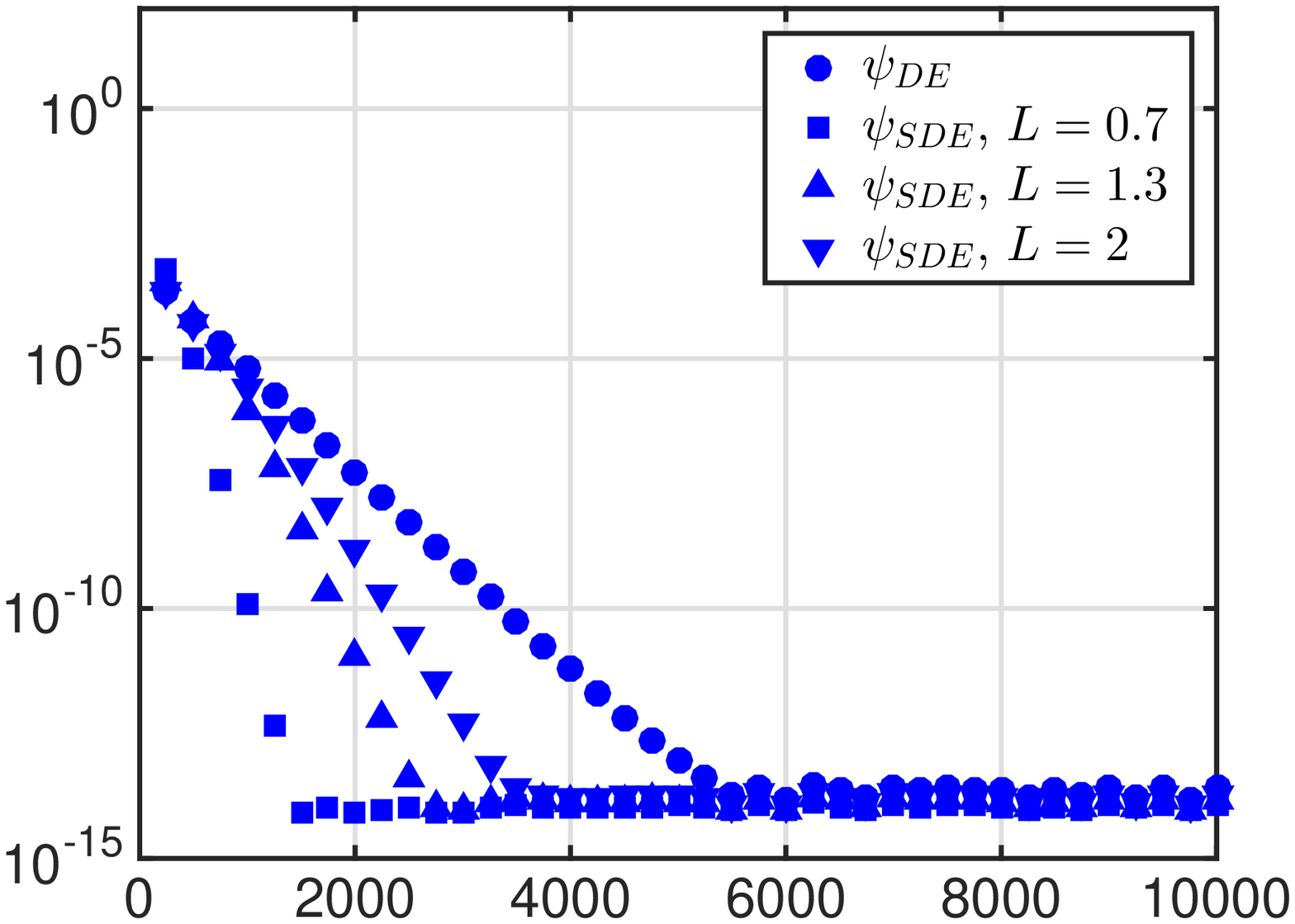}
&
\includegraphics[width=4cm]{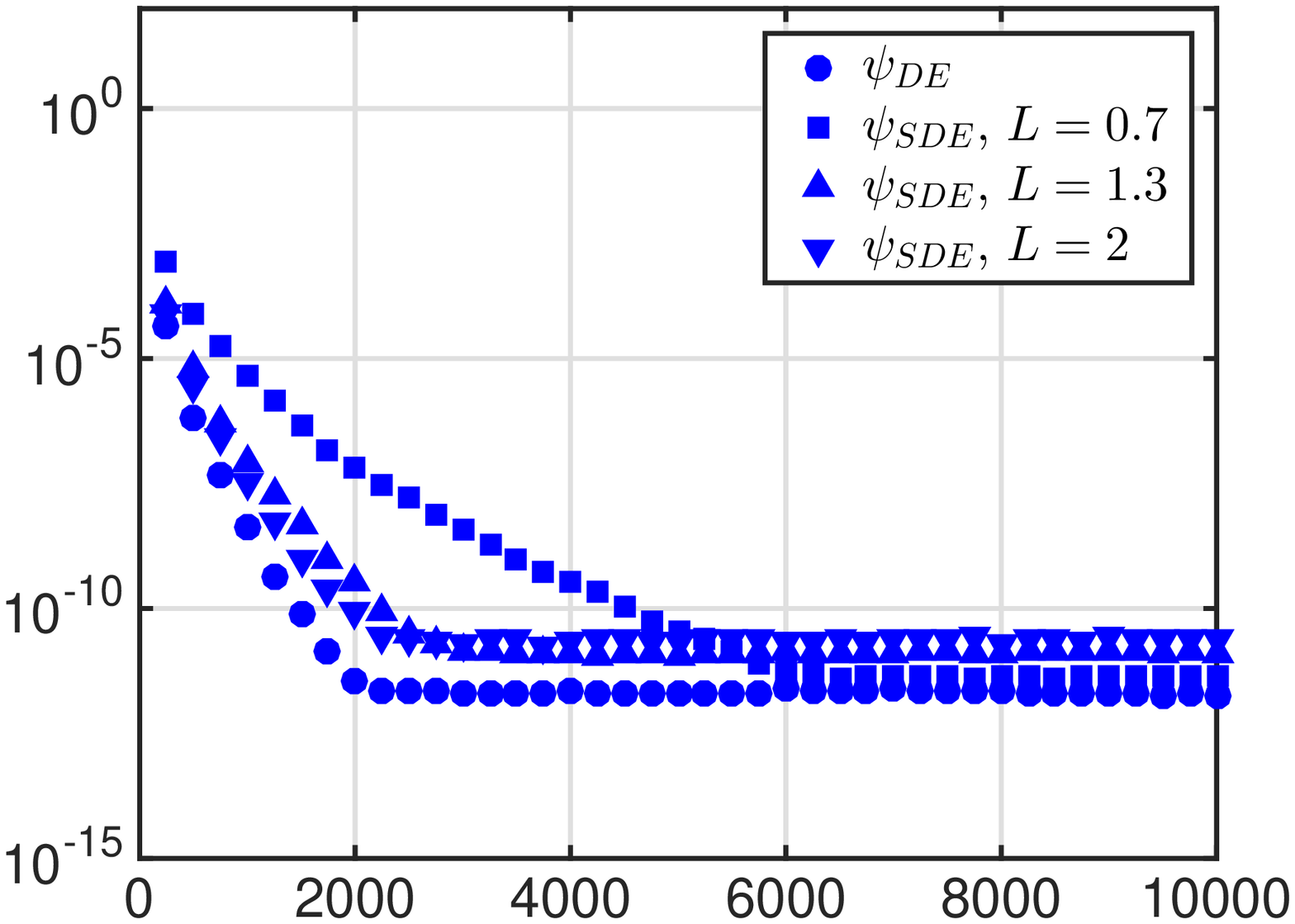}
\end{tabular}
\end{center}
\caption{ The error against $n$ for the exponential (top row) and double exponential (bottom row) maps applied to the functions $f_1(x) = x^{1/5} \exp(-800 \pi \I x)$ (left), $f_2(x) = \frac{  \sqrt{x} }  {1+100^2 (x-0.5)^2 } $ (middle) and $f_3(x) = x^{1/2} (1-x)^{3/4}  \mathrm{sn} \left(    \log \left( x^{ 5 } / ( 1-x)^{3}   \right), 1/\sqrt{2}  \right)$ (right).  Here $\mathrm{sn}(\cdot,\cdot)$ denotes the Jacobi elliptic function.  The parameters used were $L=c\sqrt{n}$, $\alpha = \alpha_0 / \sqrt{n}$, $L = W(cn)$ and $\alpha = L_0 \pi / W(cn)$ for $\psi_{E}$, $\psi_{SE}$, $\psi_{DE}$ and $\psi_{SDE}$ respectively, where the constants $c$ and $\alpha_0$ were numerically optimized.
}
\label{f:compDoubleExp}
\end{figure}

We conclude this paper with a numerical comparison of the four maps.  This is shown in Fig.\ \ref{f:compDoubleExp}.  In order to ensure a fair comparison of the various maps, the constants $c$ and $\alpha_0$ appearing in the parameter choices were numerically optimized.  This was done by varying such quantities over an appropriate range and finding the value which minimized the error for each particular function.  To make the computations feasible in a reasonable time, we have not optimized over the parameter $L$ in the parametrized maps.  Instead, we merely fix several different values of $L$.

Fig. \ref{f:compDoubleExp} compares the performance of the four maps for three different yet challenging functions to approximate.  The first is a singular oscillatory function and the second is a singular version of the classical Runge function.  The third function features both singularities and nonuniform oscillations. It is similar to a function used in \cite{Sugihara04Sinc}, yet we have changed several of the parameters to make the function more challenging to approximate.

As is evident from this figure, for the first two functions the new maps $\psi_{SE}$ and $\psi_{SDE}$ offer superior performance over the standard exponential and double-exponential maps.  In both cases the optimal choice of $L$ is the smallest, i.e.\ $L=0.7$.  For $f_1$ this is due to the results on resolution proved in this paper.  The behaviour of the error for $f_2$ is similar, since this function also grows rapidly on the shifted imaginary axis $1/2 + \I \bbR$, much like an oscillatory function.  On the other hand, this choice of $L$ leads to worse performance for $f_3$, which suggests that the convergence rate for this function is limited primarily by the singularities rather than the oscillations.  Note that for $f_3$ the new maps do not convey any advantage over the existing maps, but the performance is at least similar for the values $L=1.3$ and $L=2$.

\rem{
Implementation of the parametrized maps in finite-precision arithmetic requires a little care.  As discussed in \cite{AdcockRichardsonMappings} for $\psi_{SE}$, naive implementations of the inverse maps $\psi^{-1}_{SE}$ and $\psi^{-1}_{SDE}$ may result in cancellation errors.  Fortunately, these effects can be avoided by implementing terms such as $\exp(x) -1$ and $\log(1+x)$ using Matlab's \texttt{expm1} and \texttt{log1p} functions.  Another issue is the practical limitation on the parameter values due to the possibility for overflow or underflow.  Inspecting \R{E:PsiSE}, we see that this will generally occur for the map $\psi_{SE}$ when $\exp(\pi/\alpha)$ exceeds the largest floating point number ($\approx 10^{308}$), or in other words, $\alpha < 0.0044$.  Similarly, the same issue can occur for $\psi_{SDE}$ when $\exp(\exp(\pi/(2 \alpha)))$ exceeds $10^{308}$, or in other words, when $\alpha < 0.24$.  We have used these guidelines throughout the paper when choosing the parameter values.  Fortunately, as seen in Fig.\ \ref{f:compDoubleExp}, these barriers do not appear to hamper the performance of either map in practice, even for large values of $n$.
}

\section*{Acknowledgements}

BA acknowledges support from the Alfred P. Sloan Foundation and the Natural Sciences and Engineering Research Council of Canada through the grant 611675.
JM acknowledges support from Spanish Ministry of Econom\'ia y Competitividad through the grant TIN2014-55325-C2-2-R; and also from Spanish Ministry
of Educaci\'on, Cultura y Deporte, thanks to the mobility grant Jos\'e Castillejo CAS14/00103, JM spent three months at Simon Fraser University where this work was partially done.

\appendix

\section{Derivation of $\psi_{SDE}$}
\label{supp1}
As outlined in the main text, the map $\psiSDEi$ is formed by the composition $\psiSEi \circ g^{-1}_{}$, where $g^{-1}_{}: (-\infty, \infty)\mapsto (-\infty, \infty)$  is to be determined, and $\psiSEi$ is given by (2.11). As before, $\psiSEi$ is parameterised by a ``strip-width'' parameter $\alpha$, and we shall see that the function $g^{-1}_{}$ will be similarly parameterised by a positive real number $\gamma$.  In the paper we fix $ \gamma = \alpha$, but for clarity in the following derivation we allow $\gamma$ to be distinct from $\alpha$.
For going back and forth between the $x$ and $s$ variables, we therefore have
$
x = \psiSDEi(s; \alpha, \gamma) := \psiSEi(g^{-1}_{}(s; \gamma); \alpha),
$
and
$
s = \psiSDE(x; \alpha, \gamma) := g^{}_{}(\psiSE(x; \alpha); \gamma).
$

\begin{figure}
\begin{center}
\begin{overpic}[scale=0.42]{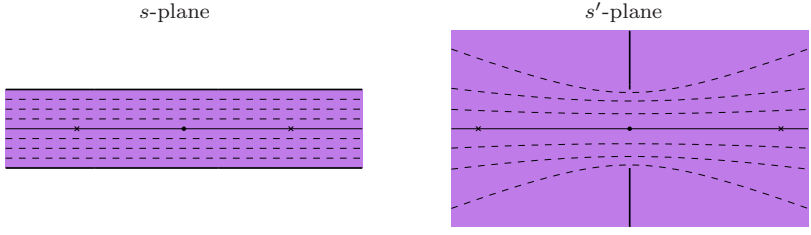}
\put(20,25){\footnotesize $s$-plane}
\put(70,25){\footnotesize $s'$-plane}
\end{overpic}
\caption{\label{fig:archsinhActionOnStrip} \small The function $\sinh$ maps the open infinite strip of half-width $\pi/2$ onto the two-slit plane $\complex \backslash \{(-\infty,-i] \, \cup \, [i,\infty)\}$.}
\end{center}
\end{figure}

To derive $g^{-1}_{}: (-\infty, \infty) \mapsto (-\infty, \infty)$, we note first that this transformation plays an analogous role to that of $\sinh$  in the standard unparameterised double-exponential map $\psiDEi(s) = \exp(\pi\sinh(s))/(1+\exp(\pi\sinh(s)))$. It is therefore reasonable to begin exploring the question of what might be an appropriate form for a resolution-optimal analogue of $\sinh$ by considering the action of $\sinh$ on an infinite strip in the complex $s$-plane. As illustrated in Figure~\ref{fig:archsinhActionOnStrip}, its action is to ``unfold'' the strip $S^{}_{\pi/2}$. As the contour lines show, this has the undesirable effect of warping functions in the complex-plane in such a way that ``strip-behaviour'' is not preserved in moving from one domain to the next, particularly in the region of the complex-plane local to $[-1/2, 1/2]$. This effect contributes to the suboptimal resolution power of $\psiDE$.

\begin{figure}
\begin{center}
\begin{overpic}[scale=0.42]{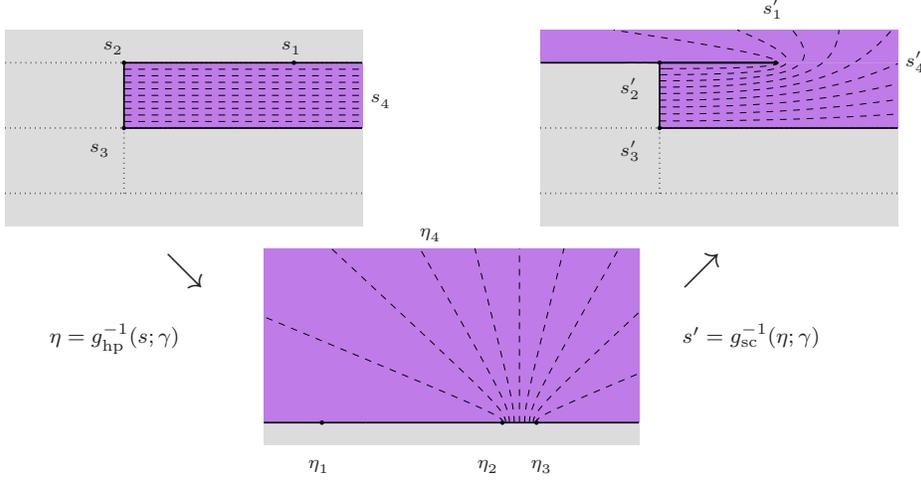}
 \put(18,21) {\Large $\searrow$}
 \put(5,14) {\footnotesize $\eta = g^{-1}_\text{hp}(s; \gamma)$} 
 \put(76,21) {\Large $\nearrow$}
 \put(76,14) {\footnotesize $s' = g^{-1}_\text{sc}(\eta; \gamma)$}
 \put(31,47) {\scriptsize $s^{}_1$} 
 \put(11,47) {\scriptsize $s^{}_2$} 
 \put(9.5,35.5) {\scriptsize $s^{}_3$} 
 \put(41,41) {\scriptsize $s^{}_4$} 
 \put(34,0) {\scriptsize $\eta^{}_1$} 
 \put(53,0) {\scriptsize $\eta^{}_2$} 
 \put(59,0) {\scriptsize $\eta^{}_3$} 
 \put(46.5,26) {\scriptsize $\eta^{}_4$} 
 \put(85,51) {\scriptsize $s'_1$} 
 \put(69,42) {\scriptsize $s'_2$} 
 \put(69,35) {\scriptsize $s'_3$} 
 \put(101,45) {\scriptsize $s'_4$} 
\end{overpic}
\caption{\label{fig:QuarterStripToQuarterFoldedStrip} \small Construction of the map $g^{-1}_{}(s;\gamma)$, which transforms a quarter-infinite strip to an ``unfolded'' quarter-infinite strip. It is derived by first mapping to the upper-half plane, before using the Schwarz-Christoffel transformation to map the destination region.}
\end{center}
\end{figure}

Our remedy is to derive a map which retains a certain amount of strip-behaviour around the real line -- an approach very much analogous to that employed in \cite{AdcockRichardsonMappings}. We proceed as in this case by constructing the transformation
\begin{equation}\label{E:g_composed}
g^{-1}_{}(s;\gamma) = g^{-1}_\text{sc}(g^{-1}_\text{hp}(s;\gamma);\gamma),
\end{equation}
where $g^{-1}_\text{hp}(s; \gamma) \ee \sinh^{2}_{} ( \pi s / 2 \gamma)$ first maps the open quarter-infinite strip
$
\{ s \s \in \s \complex: 0 \leq \Im \, s < \gamma\, , \, \Re \, s \geq 0 \}
$
to the upper-half plane, after which the map $g^{-1}_\text{sc}(\,\cdot\,;\gamma)$, takes us from the upper half-plane to an ``unfolded'' quarter-strip; see Figure~\ref{fig:QuarterStripToQuarterFoldedStrip}.

The map $g^{-1}_\text{sc}$ can be derived using the Schwarz-Christoffel approach \cite{DriscollTrefethenSchwarzChristoffel}. We begin by setting $s^{}_1 = 1/2 + \gamma i$, $s^{}_{2} = \gamma i$, $s^{}_3 = 0$, $s^{}_4 = \infty$, from which we proceed to identify the Schwarz-Christoffel prevertices $\eta_j$, vertices $s'_j$, and angles $\delta_j$ as
\eas{
\begin{array}{lll}
\eta^{}_1 = g^{-1}_\text{hp}(s^{}_{1}; \gamma),   & s'_1 = - ,     &  \delta_1 = 2, \\
\eta^{}_2 = -1,                    & s'_2 = i \gamma,          &  \delta_2 = 1/2, \\
\eta^{}_3 = 0,                    & s'_3 = 0,     &  \delta_3 = 1/2, \\
\eta^{}_4 = \infty,               & s'_4 = \infty,     &  \delta_4 = -1.
\end{array}
}
Note that we are only free to choose the location of two of the non-infinite vertices and thus we enforce the positions of $s'_2$ and $s'_3$, but not $s'_1$. Conveniently, though, the final map $g^{-1}_{}$ will have the desirable property that $\lim^{}_{\gamma \to 0} s'_1 = s^{}_1$.

\begin{figure}[ht]
\begin{center}
\begin{overpic}[scale=0.42]{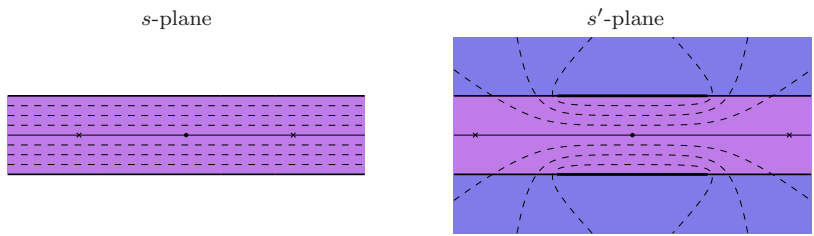}
\put(20,25){\footnotesize $s$-plane}
\put(70,25){\footnotesize $s'$-plane}
\end{overpic}
\caption{\label{fig:StripToFoldedStrip} \small For any $\gamma > 0$, $g^{-1}_{}(\,\cdot\,;\gamma)$ maps the open infinite strip $S{}_{\gamma}$ onto a ``doubly-unfolded strip''. Note the key property that ``strip-like'' behaviour in the region of the complex-plane local to the interval $[-1/2,1/2]$ is preserved in moving from one domain to the other.}
\end{center}
\end{figure}

The Schwarz-Christoffel integral is
\begin{align*}
 g^{-1}_\text{sc}(\eta; \gamma)
 = B + C \int^{\eta} (\xi - \eta^{}_1)^{-1}(\xi + 1)^{-1} \xi^{-1/2} {\rm d} \xi,
\end{align*}
for constants $B$ and $C$, which may be determined by evaluating the integral exactly, applying the composition given by \eqref{E:g_composed}, and then enforcing the two conditions $s'_2 = g(s^{}_2; \gamma)$, $s'_3 = g(s^{}_3; \gamma)$. This gives us the final map
\be{
\label{E:gfunction2}
g^{-1} (s; \gamma) = s + \frac{\gamma}{\pi} \frac{ \sinh(\pi s/ \gamma) }  { \cosh( \pi / 2 \gamma )}.
}
Note that though $g^{-1}_{}(\,\cdot\,;\gamma)$ was constructed using considerations on a quarter-strip, by the Schwarz reflection principle the map can be extended across the boundaries such that it is in fact valid across the entire strip. This can be seen in Figure~\ref{fig:StripToFoldedStrip}.

\section{Analysis of the map $\psi_{SDE}$}
\label{supp2}
Here we present some technical lemmas in the analysis of the map $\psi^{-1}_{SDE}$.

\lem{
\label{l:SDEanalytic}
The map $\psi^{-1}_{SDE}$ is analytic in the domain $S_{\beta}$ for any $\beta < \alpha/2$.
}
\prf{
Since $\psi^{-1}_{SDE} = \psi^{-1}_{SE} \circ g^{-1}$ and $g^{-1}$ is entire it suffices to consider $\psi^{-1}_{SE}$.  Note that $\psi^{-1}_{SE}$ fails to be analytic if and only if one of the following holds:
\begin{itemize}
\item[(i)] $1 + \exp(\pi(s-1/2) / \alpha) = 0$,
\item[(ii)] $\Im (z) = 0$ and $\Re(z) \leq 0$, where $z = \frac{1+\exp(\pi(s+1/2)/\alpha)}{1+\exp(\pi(s-1/2)/\alpha)}$.
\end{itemize}
Condition (i) holds if and only if $s=1/2+(2k+1) \alpha \I$, $k \in \bbZ$.  Now consider condition (ii).  If $s = x + \I y$ then
\eas{
\Re(z) &= \frac{\exp(\pi/(2\alpha)) \left ( \cos(\pi y / \alpha) \cosh(\pi/(2\alpha)) + \cosh(\pi x / \alpha) \right )}{\cos(\pi y / \alpha) + \cosh \pi(x-1/2)/\alpha) },
\\
\Im(z) &= \frac{(\exp(\pi/\alpha)-1) \sin(\pi y / \alpha)}{2\left ( \cos(\pi y / \alpha) + \cosh ( \pi(x-1/2) / \alpha) \right )}.
}
Hence $\Im(z) = 0$ if and only if
\bes{
y = k \alpha,\qquad k \in \bbZ,
}
in which case we have
\bes{
\Re(z) = \frac{\exp(\pi / (2\alpha)) \left ( (-1)^k \cosh(\pi / (2\alpha)) + \cosh(\pi x / \alpha) \right )}{(-1)^k + \cosh ( \pi ( x-1/2 ) / \alpha)}
}
When $k$ is even $\Re(z) |_{y = k \alpha}$ is always positive.  If $k$ is odd, then $\Re(z) |_{y = k \alpha}$ is negative if and only if $|x| < 1/2$.  Hence condition (ii) holds if and only if
\bes{
y = (2k+1) \alpha,\quad |x| < 1/2,\qquad k \in \bbZ.
}
Combining conditions (i) and (ii) together we now deduce that $\psi^{-1}_{SE}$ fails to be analytic at a point $s = x+\I y$ if and only if
\bes{
|x| \leq 1/2,\qquad y = (2k+1)\alpha,\ k \in \bbZ.
}
From this, we see that $\psi^{-1}_{SDE}$ fails to be analytic at a point $z = x+\I y \in S_{\beta}$ if and only if
\be{
\label{bad_cond}
\left | \Re (g^{-1}(x+\I y) )\right | \leq 1/2,\qquad \Im (g^{-1}(x+\I y)) = (2k+1) \alpha,\qquad k \in \bbZ.
}
By explicit calculation
\eas{
\Re(g^{-1}(x+\I y)) &= x + \frac{\alpha \cos(\pi y / \alpha) \sinh(\pi x/ \alpha)}{\pi \cosh(\pi / (2\alpha))},
\\
\Im(g^{-1}(x+\I y)) & = y + \frac{\alpha \sin(\pi y / \alpha) \cosh(\pi x/ \alpha)}{\pi \cosh(\pi / (2\alpha))}.
}
Observe that \R{bad_cond} cannot hold for $y=0$.  If $y \neq 0$, then it follows that \R{bad_cond} holds only if
\be{
\label{bad_cond_x}
\cosh(\pi x / \alpha) = \frac{\pi \cosh(\pi / (2\alpha))}{\alpha \sin ( \pi y / \alpha)} \left ( (2k+1) \alpha - y \right ).
}
Without loss of generality, let $k \geq 0$.  Since $x + \I y \in S_{\beta}$, and therefore $|y| < \beta < \alpha /2< (1-1 /\pi) \alpha$ we see that \R{bad_cond_x} only possibly holds for $y > 0$.  Moreover, we must have
\bes{
\cosh(\pi x / \alpha) > \cosh(\pi / (2\alpha)),
}
and therefore $|x| > 1/2$.  But if $0 < y < \alpha / 2$ and $|x| > 1/2$ then $| \Re(g^{-1}(x+\I y) | > 1/2$.  Hence \R{bad_cond} cannot hold, as required.
}

\lem{
\label{PsiSDEres_step}
For $\alpha > 0$ let $\psi_{SDE}$ be the mapping given by (2.12).  Suppose that $L$ is uniformly bounded above and $L-1/2>0$ is uniformly bounded away from zero for all $\alpha$.  Then, for the function $f(x) = \exp(-2 \pi \I \omega x)$ we have
\bes{
M_{\psi}(\sigma \alpha ; f ) \leq \exp(2 \pi \omega\sigma ( \alpha + o(\alpha) ) ),\qquad \alpha \rightarrow 0,
}
and
\bes{
N_\psi(\sigma \alpha , L ; f) \leq 2 \pi \omega ( 1 + o(1) ),\qquad \alpha \rightarrow 0,
}
uniformly in $\omega \geq 0$ for $0 < \sigma \leq \sigma^*$, where $\sigma^* \approx 0.265$ is the unique root of $\sigma \pi + \sin(\sigma \pi) = \pi/2$ in $0 < \sigma < 1/2$.
}
\prf{
We first consider $M_\psi$:
\bes{
M_{\psi} = \exp \left ( 2 \pi \omega \sup_{x \in \bbR} \Im \psi^{-1}(x \pm \I \sigma \alpha) \right ) = \exp \left ( 2 \pi \omega \sup_{x \geq 0}\Im \psi^{-1}(x \pm \I \sigma \alpha) \right ).
}
where in the second step we use the symmetry relation (2.2).  As $\alpha \rightarrow 0$, note that
\bes{
g^{-1}(x\pm\I \sigma \alpha) \sim x \pm \I \sigma \alpha + \frac{2\alpha}{\pi} \exp(-\pi/(2\alpha)) \sinh(\pi x / \alpha \pm \I \sigma \pi),
}
uniformly in $x \geq 0$ and $0 < \sigma < 1/2$.  By definition
\bes{
\Im \psi^{-1}(x \pm \I \sigma \alpha) = \frac{\alpha}{\pi} \Im \log \left ( \frac{1+\exp(\pi (g^{-1}(x\pm\I \sigma \alpha)  +1/2)/\alpha)}{1+\exp(\pi (g^{-1}(x\pm\I \sigma \alpha)  -1/2)/\alpha)} \right ).
}
For the numerator, we have
\bes{
\Im \log \left ( 1+\exp(\pi (g^{-1}(x\pm\I \sigma \alpha)  +1/2)/\alpha) \right ) \sim \pm \theta(x,\alpha , \sigma),
}
where
\bes{
\theta(x,\alpha,\sigma) = \sigma \pi + 2 \exp(-\pi / (2\alpha)) \cosh(\pi x / \alpha) \sin(\sigma \pi).
}
For the denominator, we have
\bes{
1+\exp(\pi (g^{-1}(x\pm\I \sigma \alpha)  -1/2)/\alpha) \sim 1 + R(x,\alpha,\sigma) \exp( \pm \I \theta(x,\alpha,\sigma)),
}
where
\eas{
R(x,\alpha,\sigma) & = \exp \left ( \pi (x-1/2) / \alpha + 2 \exp(-\pi / (2 \alpha)) \sinh(\pi x / \alpha) \cos(\sigma \pi) \right ).
}
Hence, after dividing,
\eas{
\Im \psi^{-1}(x\pm \I \sigma \alpha) &\sim - \frac{\alpha}{\pi} \Im \log \left [ R(x,\alpha,\sigma) + \exp(\mp \I \theta(x,\alpha,\sigma)) \right ],
}
and therefore
\bes{
\Im \psi^{-1}(x\pm \I \sigma \alpha)  \leq \frac{\alpha}{\pi} \arctan \left | \frac{\sin(\theta(x,\alpha,\sigma))}{R(x,\alpha,\sigma) + \cos(\theta(x,\alpha,\sigma))} \right |.
}
We now consider two cases:

\vspace{1pc} \noindent Case 1 ($x \geq 1/2$): As $\alpha \rightarrow 0$, we have
\eas{
R(x,\alpha,\sigma) &\sim \exp \left ( \pi(x-1/2) / \alpha + \exp(\pi(x-1/2) / \alpha) \cos(\sigma \pi) \right )
\\
\theta(x,\alpha,\sigma) &\sim \sigma \pi + \exp(\pi(x-1/2)/\alpha) \sin (\sigma \pi),
}
uniformly in $x \geq 1/2$.  Hence if $y = \pi(x-1/2)/\alpha \geq 0$, we deduce that
\bes{
\frac{\sin(\theta(x,\alpha,\sigma))}{R(x,\alpha,\sigma) + \cos(\theta(x,\alpha,\sigma))} \sim \frac{\sin(\sigma \pi + \exp(y) \sin(\sigma \pi))}{\exp(y + \exp(y) \cos(\sigma \pi)) + \cos(\sigma \pi + \exp(y) \sin(\sigma \pi))}.
}
It now follows from Lemma \ref{l:technical_step} that
\bes{
\sup_{x \geq 1/2} \Im \psi^{-1}(x\pm \I \sigma \alpha) \leq \sigma \alpha + o(\alpha),\qquad \alpha \rightarrow 0.
}

\vspace{1pc} \noindent Case 2 ($0 \leq x < 1/2$): Note first that
\bes{
\theta(x,\alpha,\sigma) \sim \sigma \pi + \exp(\pi(x-1/2) /\alpha) \sin(\sigma \pi),\qquad \alpha \rightarrow 0,
}
uniformly in $0 \leq x < 1/2$.  Therefore
\bes{
0 < \sigma \pi < \theta(x,\alpha,\sigma) \leq \sigma \pi + \sin(\sigma \pi),
}
for all small $\alpha$ and all $0 \leq x < 1/2$.  In particular, $\sigma \pi \leq \theta(x,\alpha,\sigma) \leq \pi/2$ provided $0 < \sigma < \sigma^*$.  For such values of $\sigma$, it follows that
\bes{
\frac{\sin(\theta(x,\alpha,\sigma))}{R(x,\alpha,\sigma) + \cos(\theta(x,\alpha,\sigma))} > 0.
}
With this in hand, we now claim that
\be{
\label{PsiSDEres_claim}
\frac{\sin(\theta(x,\alpha,\sigma))}{R(x,\alpha,\sigma) + \cos(\theta(x,\alpha,\sigma))} \leq \tan(\sigma \pi),\qquad \alpha \rightarrow 0,
}
uniformly in $0 \leq x < 1/2$ and $0 < \sigma \leq \sigma^*$.  We have
\eas{
\frac{\sin(\theta(x,\alpha,\sigma))}{R(x,\alpha,\sigma) + \cos(\theta(x,\alpha,\sigma))} = \left ( \frac{\frac{\sin\left(\theta(x,\alpha,\sigma) - \sigma \pi \right )}{\sin(\sigma \pi)} + \cos(\theta(x,\alpha,\sigma))}{R(x,\alpha,\sigma) + \cos(\theta(x,\alpha,\sigma))} \right ) \tan(\sigma \pi),
}
and, since the denominator is positive it suffices to show that
\bes{
\frac{\sin\left(\theta(x,\alpha,\sigma) - \sigma \pi \right )}{\sin(\sigma \pi)} \leq R(x,\alpha,\sigma),\qquad \alpha \rightarrow 0.
}
By definition and the fact that $ \theta(x,\alpha,\sigma) > \sigma \pi$, we see that
\eas{
\frac{\sin\left(\theta(x,\alpha,\sigma) - \sigma \pi \right )}{\sin(\sigma \pi)} &\leq 2 \exp(-\pi / (2 \alpha)) \cosh(\pi x / \alpha )
\\
& \sim \exp(\pi(x-1/2)/\alpha)
\\
& \leq R(x,\alpha,\sigma),
}
as required.  Thus \R{PsiSDEres_claim} is proved, and it now follows that
\bes{
\sup_{0 \leq x \leq 1/2} \Im \psi^{-1}(x\pm\I \sigma \alpha) \leq \sigma \alpha + o(\alpha),\qquad \alpha \rightarrow 0.
}
Combining this with the previous case, we now finally arrive at
\bes{
\sup_{x \geq 0} \Im \psi^{-1}(x\pm\I \sigma \alpha) \leq \sigma \alpha + o(\alpha),\qquad \alpha \rightarrow 0,
}
and this completes the proof for $M_\psi$.

We now consider $N_{\psi}$.  Since $\tau = 1$, we have
\bes{
N_{\psi} \leq 2 \pi \omega  \exp \left ( 2 \pi \omega \sup_{0 \leq \theta \leq 2 \pi} \psi^{-1}(L+\sigma \alpha \E^{\I \theta}) \right ).
}
Since $\Re(L + \sigma \alpha \E^{\I \theta}) > 1/2$ for all small $\alpha$, we may argue as above and deduce that
\bes{
\Im \psi^{-1}(L+\sigma \alpha \E^{\I \theta}) \sim 0,\qquad \alpha \rightarrow 0.
}
This now gives the result.
}

\lem{
\label{l:technical_step}
Let $G(y,\sigma) = \frac{\sin(\sigma \pi + \exp(y) \sin(\sigma \pi))}{\exp ( y + \exp(y) \cos(\sigma \pi)) + \cos(\sigma \pi + \exp(y) \sin (\sigma \pi))}$.  Then
\bes{
\sup_{y \geq 0 }\left | G(y,\sigma) \right | \leq \tan(\sigma \pi),\qquad \forall 0 < \sigma < 1/2.
}
}
\prf{
We first show that $G(y,\sigma) \leq \tan(\sigma \pi)$.   Fix $0 < \sigma < 1/2$.  Rearranging and simplifying, we see that this is equivalent to
\bes{
g(y) = \sin(\exp(y) \sin(\sigma \pi))/ \sin(\sigma \pi)  \leq h(y) =  \exp ( y + \exp(y) \cos(\sigma \pi)).
}
Observe that $g(0) = \sin(\sin(\sigma \pi))/\sin(\sigma \pi) < 1$ and $h(0) = \exp(\cos(\sigma \pi)) > 1$.  Also,
\bes{
g'(y) = \exp(y) \cos(y \sin(\sigma \pi)) \leq \exp(y),
}
whereas
\bes{
h'(y) = \left ( 1 + \exp(y) \cos(\sigma \pi) \right ) \exp(y + \exp(y) \cos(\sigma \pi)) \geq \exp(y).
}
Hence $g(y) \leq h(y)$, $\forall y \geq 0$, as required.

We now show that $G(y,\sigma) \geq -\tan(\sigma \pi)$.  Rearranging and simplifying once more, this is equivalent to
\bes{
g(y) = \sin(2\sigma \pi + \exp(y) \sin(\sigma \pi)) / \sin(\sigma \pi) \geq h(y) = - \exp(y + \exp(y) \cos(\sigma \pi)).
}
Note that $g(0) = \sin(\sigma \pi + \sin(\sigma \pi))/\sin(\sigma \pi) > 0$.  Conversely,  $h(0) = - \exp(\cos(\sigma \pi)) < 0$.  Also,
\bes{
g'(y) = \exp(y) \cos(2 \sigma \pi + \exp(y) \sin(\sigma \pi)) \geq - \exp(y)
}
and
\bes{
h'(y) = - \left ( 1 + \exp(y) \cos(\sigma \pi) \right ) \exp( y + \exp(y) \cos(\sigma \pi) ) \leq -\exp(y).
}
Thus $g(y) \geq h(y)$, $\forall y \geq 0$, as required.
}

\bibliographystyle{siamplain}
\small
\bibliography{FourMapsSingul}

\end{document}